%% file: Metric_SG.tex
\documentclass[11pt, dvipdfmx]{amsart}
\usepackage{amsmath, amsthm, amscd, amssymb, enumerate, verbatim, graphicx}
\usepackage[all]{xy}
\usepackage{xcolor}
\usepackage[T1]{fontenc}

\input macro.tex

\begin{document}

\title[]{Note on the induced metric on power sets of groups} 

\author[Tatsuhiko Yagasaki]{Tatsuhiko Yagasaki}
\address{Professor Emeritus \\
Kyoto Institute of Technology \\
Kyoto, Japan}
\email{yagasaki@kit.ac.jp}


\begin{abstract}{}
In this note we clarify general properties of the Hausdorff-like metric on the power sets ${\cal S}(G)$ of a group $G$ induced from word length norm 
and obtain some results on quasi-isometries between some subspaces of ${\cal S}(G)$ and ${\cal S}(H)$ for 
a group epimorphism $f : G \to H$, when ${\cal S}({\rm Ker}\,f)$ is metrically bounded. 
\end{abstract} 

\maketitle

\thispagestyle{empty}

\tableofcontents

\baselineskip 6mm 

\section{Introduction}

The study of boundedness of groups aims a quantitative understanding of various norms and metrics on groups or   
some spaces associated to them (cf. \cite{BIP, Fuj, KLM1, KKKMMO, Ts2, Ts3, Ts4}, etc). 
In \cite{KLM1, KLM2} J.~ K\k{e}dra, A.~Libman and B.~Martin  %Jarek K\k{e}dra, Assaf Libman and Ben Martin 
studied strong and uniform boundedness of finitely normally generated groups based on the word lengths $\| \ \|_S$ with respect to finite normally generating sets $S$ on a group. 
For a study of unbounded groups, in \cite{KKKMMO} they consider coarse type of Tsuboi metric spaces associated to relatively simple groups and also  
after \cite{KLM1, KLM2} Morimichi Kawasaki have studied the metric on the set of finite normally generating sets of a finitely normally generated group modulo symmetrization and conjugation, obtained by symmetrizing the word length norm in \cite{KLM1, KLM2}. 

Our goal in this subject is a complehensive study of boundedness of bundle diffeomorphism groups and diffeomorphism groups of manifold pairs. 
In \cite{FY1, FY2} we have studied the cases of bundle diffeomorphism groups over a circle  and diffeomorphism groups of manifold pairs, in which   
the submanifold is a finite disjoint union of circles, and are in preparation of successive papers in the higher dimensional cases. 
As a preliminary of these works, in this note we try to make a framework and a foundation under a general setting and clarify general properties  on the extended metric on the power set ${\cal S}(G)$ of a group $G$,
which is induced from the word lengths with respect to subsets of $G$. 
Here, the power set ${\cal S}(G)$ is the set of all subsets of $G$ and 
we treat the wide class of extended (additive/multiplicative) assymmetric metrics with values in $\IR_{\geq 0} \cup \{ \infty \}$. 
Each $S \in {\cal S}(G)$ induces the word length $\nu_S$ on $G$ and its Hausdorff version $(\nu_S)_H$ on ${\cal S}(G)$, 
from which obtained are the multiplicative/additive assymmetric metrics $\nu(S, T) = (\nu_S)_H(T)$ and $\rho(S,T) = \log \nu_H(S,T)$ 
(cf. \cite[Section 2]{KLM1}). 
The usual symmetric metrics $\widehat{\nu}$ and $\widehat{\rho}$ are readily obtained by symmetrizing $\nu$ and $\rho$. 
We list formal properties of these extended metrics. 
In this general setting, one of main problems is to understand relations as metric spaces among (some subspaces of) the power sets ${\cal S}(G)$, ${\cal S}(H)$ and ${\cal S}(K)$ for a group epimorphism $f : G \to H$ with $K = {\rm Ker}\,f$. 
We obtain some weak results on quasi-isometries between some subspaces of ${\cal S}(G)$ and ${\cal S}(H)$ when ${\cal S}(K)$ is metrically bounded. 

Although the materials in this note are seemed to be essentially known for the authors in \cite{KLM1, KLM2},  
we expect that it is useful for our works to clarify general properties of the metrics in ${\cal S}(G)$ in this general setting. 
As an advantage of this framework
we can treat various subsets of a group simultaneously within the ambient extended metric space. 
We also extend this framework to the cases of group actions and sets with binary operations. 

\section{Baiscs on asymmetric metrics and Hausdorff metrics}

This section is devoted to generalities on asymmetric metrics and Hausdorff metrics on powersets. 

\subsection{Asymmetric metrics} \mbox{} 

Let $X$ be a set.  

\bdefn\label{def_metric}
For a function $d : X \times X \lra \IR_{\geq 0} \cup \{ \infty \}$ consider the following conditions : 
\bit 
\itemI $d(x,x) = 0$ \ \ $(\upfa x \in X)$ 
\hspp (i)$'$ \ [\,$d(x,y) = d(y,x) = 0$ \LLRA $x = y$\,] \ \ $(\upfa x,y \in X)$ 
\itemII $d(x,y) = d(y,x)$ \ \ $(\upfa x,y \in X)$
\itemiii $d(x,z) \leq d(x,y) + d(y,z)$ \ \ $(\upfa x,y,z \in X)$  
\eit 
\hsf We say that the function $d$ is 
\benum 
\item[1)] an asymmetric metric on $X$ if $d$ satisfies (i) and (iii), 
\item[2)] nondegenerate if it satisfies (i)$'$, 
\item[3)] a psuedo metric on $X$ if it satisfies (i), (ii) and (iii), 
\item[4)] a metric on $X$ if it satisfies (i)$'$, (ii) and (iii). 
\eenum 
\edefn  

\bdefn\label{def_metric}
For a function $\rho : X \times X \lra \IR_{\geq 1} \cup \{ \infty \}$ consider the following conditions : 
\bit 
\itemI $\rho(x,x) = 1$ \ \ $(\upfa x \in X)$ 
\hspp (i)$'$ \ [\,$\rho(x,y) = \rho(y,x) = 1$ \LLRA $x = y$\,] \ \ $(\upfa x,y \in X)$ 
\itemII $\rho(x,y) = \rho(y,x)$ \ \ $(\upfa x,y \in X)$
\itemiii $\rho(x,z) \leq \rho(x,y) \cdot \rho(y,z)$ \ \ $(\upfa x,y,z \in X)$
\eit 
\hsf We say that the function $\rho$ is  
\benum 
\item[1)] a multiplicative asymmetric metric on $X$ if $d$ satisfies (i) and (iii), 
\item[2)] nondegenerate if it satisfies (i)$'$, 
\item[3)] a multiplicative psuedo metric on $X$ if it satisfies (i), (ii) and (iii), 
\item[4)] a multiplicative metric on $X$ if it satisfies (i)$'$, (ii) and (iii).  
\eenum 
Briefly, ``an asymmetric metric'' is abbreviated as ``an as metric''. 
\edefn

For notational simplicity we use the following notations. 
For any subset $J \subset \IR$ and $r \in \IR$ we set 
$J_{\geq r} := \{ x \in J \mid x \geq r \}$ and $\widetilde{J}_{\geq r} := J_{\geq r} \cup \{ \infty \}$. 
For sets $X$ and $Y$ the symbol $F(X,Y)$ denotes the set of functions from $X$ to $Y$. 
For $f \in F(X, \widetilde{J}_{\geq r})$ we write $f < \infty$ if $f$ takes its values in $J_{\geq r}$.

\bnotb 
The notations ${\cal M}(X; J)$, ${\cal M}_{p}(X; J)$ and ${\cal M}_{as}(X; J)$ denote 
the sets of metrics, psuedo metrics and as metrics on a set $X$ 
with values in $J \subset \widetilde{\IR}_{\geq 0}$, respectively. 
Similarly, ${\cal M}_m(X; I)$, ${\cal M}_{m,p}(X; I)$ and ${\cal M}_{m, as}(X; I)$ 
denote the sets of multiplicative metrics, psuedo metrics and as metrics a set $X$ 
with values in $J \subset \widetilde{\IR}_{\geq 1}$, respectively. 
The symbols $J$ and $I$ are omitted when $J = \widetilde{\IR}_{\geq 0}$ and $I = \widetilde{\IR}_{\geq 1}$. 
For $d \in {\cal M}_{as}(X)$ 
its restriction to $A \subset X$, 
$d|_{A \times A} \in  {\cal M}_{as}(A)$, is simply denoted by $d|_A$.  
\enot  

In this note we work in the set ${\cal M}_{as}(X)$. 
Since $d \in {\cal M}_{as}(X)$ may be degenerate or non-symmetric and may take the value $\infty$, 
we list so basic properties in order to remind and confirm them, while 
their proofs are omitted since those follows directly form their definitions and routine arguments. 

\bnotb Let $d \in F(X^2, \widetilde{\IR}_{\geq 0})$. 
The symmetrization of $d$ is defined by \\[1mm] 
\hsppp $\widehat{d} : X \times X \lra \widetilde{\IR}_{\geq 0}$ : \ $\widehat{d}(x,y) = \max \{ d(x,y), d(y,x) \}$. 
\vskip 1mm 
If $d \in {\cal M}_{as}(X)$, then (i) \ $\widehat{d} \in {\cal M}_{p}(X)$ and (ii) \ $\widehat{d} \in {\cal M}(X)$  if and only if $d$ is nondegenerate. 
\enot 

We recall some notions on (quasi-)isometries.

\bnotb Suppose $f : (X,d) \to (Y, \rho)$ is a map between (additive/multiplicative) as metric spaces. 
We use the following terminologies. 
\benum 
\item 
\bit 
\itemI $f$ is isometric if $\rho(f(x), f(x')) = d(x,x')$ \ $(\upfa x,x' \in X)$. \\
If $f$ is isometric and $d$ is nondegenerate, then $f$ is injective. 
\itemII $f$ is an isometry if $f$ is an isometric bijection. 
\eit 
\item $f$ is a quasi-isometric map, if it satisfies the following condition : \\
\hsp $(\dagger)$ \ $\rho(f(x), f(x')) \leq \lambda d(x,x') + C$ \ $(\upfa x,x' \in X)$ \ \ for some $\lambda \in \IR_{\geq 1}$ and $C \in \IR_{\geq 0}$.

\item[] $\circ$ When $d$ and $\rho$ are multiplicative as metrics, the condition $(\dagger)$ is replaced as follows : \\
\hsp $(\ddagger)$ \ $\rho(f(x), f(x')) \leq D d(x,x')^\lambda$ \ $(\upfa x,x' \in X)$ \ \ for some $\lambda \in \IR_{\geq 1}$ and $D \in \IR_{\geq 1}$.

\item Let $r \in \IR_{\geq 0}$. Two maps $h, h' : X \to (Y, \rho)$ are $r$-close if \ $\widehat{\rho}(h(x), h'(x)) \leq r$ \ $(\upfa x \in X)$. 
\item $f$ is a quasi-isometry if $f$ is a quasi-isometric map and there exists a quasi-isometric map \\
\hsh $g : (Y, \rho) \to (X,d)$ such that \ $gf$ is $r$-close to $\id_X$ and $fg$ is $r$-close to $\id_Y$ for some $r \in \IR_{\geq 0}$. 
\eenum 
\enot

\bnotb Let $W$ be a set and $f,g \in F(W, \widetilde{\IR}_{\geq 0})$. 
\benum 
\item $f \leq_q g$ \LLRA $f \leq \lambda g + C$ \ for some $\lambda \in \IR_{\geq 1}$ and $C \in \IR_{\geq 0})$. 
\vskip 0.5mm 
\item $f \sim_q g$ \ (quasi-equivalent) \LLRA $f \leq_q g$ \ and \ $g \leq_q f$ 
\eenum 
\enot

Note that for $d,d' \in {\cal M}(X)$, $d \sim_q d'$ \ if and only if \ $\id_X : (X,d) \lra (X, d')$ is a quasi-isometry. 

We can introduce a natural metric to the set ${\cal M}_{as}(X) \subset F(X^2, \widetilde{\IR}_{\geq 0})$. 

\bnot\label{not_F} Let $W$ be any set. 
\benum 
\item The set $F(W, \widetilde{\IR}_{\geq 0})$ has the following 
nondegenerate as metric. \\[2mm] 
\hsp $\mu(f,f') := \left\{ \hspace*{-1mm} \bary[c]{ll}
\inf \{ r \in \IR_{\geq 0} \mid f' \leq f + r \} & (\upexi  r \in \IR_{\geq 0} \ \ \mbox{s.t.} \ f' \leq f + r) \\[1.5mm]
\ \infty & (\mbox{otherwise}) 
\eary \right.$
\vskip 2mm 

\item The multiplicative version of $\mu$ is defined by \\[1mm] 
\hsp $\lambda : F(W, \widetilde{\IR}_{\geq 1})^2 \lra \widetilde{\IR}_{\geq 1}$ : \ 
$\lambda(g,g') := \left\{ \hspace*{-1mm} \bary[c]{ll}
\inf \{ s \in \IR_{\geq 1} \mid g' \leq s g \} & (\upexi  s \in \IR_{\geq 1} \ \ \mbox{s.t.} \ g' \leq s g) \\[2mm]
\ \infty & (\mbox{otherwise}) 
\eary \right.$
\vskip 3mm 

\item The additive as metric $\mu$ and the multiplicative as metric $\lambda$ correspond under the reciprocal bijections : \\[1mm] 
\hsp \hsh $\bary[t]{@{}c@{ \ }c@{ \ }c@{ \ }c@{}}
\exp_\ast : & F(W, \widetilde{\IR}_{\geq 0}) & \lra & F(W, \widetilde{\IR}_{\geq 1}) \\[0.5mm] 
& f &  & \exp f  
\eary$ \hsp 
$\bary[t]{@{}c@{ \ }c@{ \ }c@{ \ }c@{}}
\log_\ast : & F(W, \widetilde{\IR}_{\geq 1}) & \lra & F(W, \widetilde{\IR}_{\geq 0}) \\[0.5mm] 
& g &  & \log g 
\eary$ 
\vskip 1mm 

\bit 
\itemI $\mu = \log \lambda\,(\exp_\ast \times \exp_\ast)$
\vskip 1mm 
\itemII For $f, f' \in F(W, \widetilde{\IR}_{\geq 0})$ and $r \in \IR_{\geq 0}$ \\[0.5mm] 
\hsp $\mu(f,f') \leq r$ \LLRA $f' \leq f + r$ \LLRA $e^{f'} \leq e^r e^f$ \LLRA $\lambda(e^f, e^{f'}) \leq e^r$
\eit 
\eenum 
\enot

\bdefnb The set ${\cal M}_{as}(X)$ inherits the nondegenerate as metric $\mu|_{{\cal M}_{as}(X)}$ 
from $(F(X^2, \widetilde{\IR}_{\geq 0}), \mu)$. 
\bit 
\item[] \hsh (i) \ $\mu(d, d') \leq r$ \LLRA $d' \leq d + r$ \hsh $(d, d' \in {\cal M}_{as}(X), r \in \IR_{\geq 0})$ 
\eit 
\edefn 

\bnotb Suppose $(X,d)$, $(Y, \rho)$ are as metric spaces. 
\benum

\item We have an as metric $\mu \in {\cal M}_{as}(X \times Y)$ : $\mu((x,y), (x',y')) = \max \{ d(x,x'), \rho(y,y') \}$. \\
If $d$ and $\rho$ are nondegenerate, then $\mu$ is also nondegenerate. 
\eenum 
\enot 

\subsection{Hausdorff metrics on the power sets} \mbox{} 

Next we recall basic notions on Hausdorff metrics. Let $X$ be a set. 
The symbol ${\cal S}(X)$ denote the power set of $X$, that is, the set of all subsets of $X$. 
We also set ${\cal S}_f(X) := \{ S \in {\cal S}(X) \mid |S| < \infty \}$.
We regard as $X \subset {\cal S}(X)$ by identiying each point $x \in X$ with the singleton $\{ x \} \in {\cal S}(X)$. 
Each $d \in {\cal M}_{as}(X)$ determines the Hausdorff metric $d_H\in {\cal M}_{as}({\cal S}(X))$. 
 The $r$ balls and $r$ neighborhoods with respect to $d$ 
are denoted by the following notations.  
 
\bnotb Let $d \in {\cal M}_{as}(X)$, $x \in X$ and $A \subset X$.
 \benum 
\item $B_d(x, r) := \{ y \in X \mid d(x,y) \leq r \}$ \hsh $B_d(A, r) := \bigcup_{a \in A} B_d(a, r)$ \hsp $(r \in \IR_{\geq 0})$
\vskip 1mm 
\item $B_d(x, \infty) := \bigcup_{r \in \IR_{\geq 0}} B_d(x, r)$ \hspace{12mm} 
$B_d(A, \infty) := \bigcup_{a \in A} B_d(a, \infty) = \bigcup_{r \in \IR_{\geq 0}} B_d(A, r)$ 
\eenum 
We omit the symbol $d$ when it is understood from the context. 
\enot

\bfactb 
\benum 
\item $x \in B(x,0)$, \ $A \subset B(A,0)$ \hsp (2) \ 
$A \subset A'$ \LRA $B(A, r) \subset B(A', r)$ \ \ $(r \in \widetilde{\IR}_{\geq 0})$
\vskip 1mm 
\item[(3)] $B(B(A, r), s) \subset B(A, r+s)$ \ \ $(r,s \in \widetilde{\IR}_{\geq 0})$
\eenum 
\efact 

\bdefnb Each $d \in {\cal M}_{as}(X)$ induces an as metric $d_H \in {\cal M}_{as}({\cal S}(X))$ defined as follows. 
\benum 
\item $d_H : {\cal S}(X) \times {\cal S}(X) \lra \widetilde{\IR}_{\geq 0}$ : \\[2mm] 
\hsp $d_H(A, A') := \left\{ \hspace*{-1mm} \bary[c]{ll}
\inf \{ r \in \IR_{\geq 0} \mid A' \subset B_d(A, r) \} & (\upexi  r \in \IR_{\geq 0} \ \ \mbox{s.t.} \ A' \subset B_d(A, r)) \\[2mm]
\ \infty & (\mbox{otherwise}) 
\eary \right.$
\vskip 2mm 
\item $B_H(A, r) := \{ x \in X \mid d_H(A, x) \leq r \}$ \hsp $B_H(A, \infty) := \bigcup_{r \in \IR_{\geq 0}} B_H(A, r)$ 
\vskip 1mm 
\item $A \subset X$ : $d$-closed \LLRA $B_H(A, 0) = A$ \hsp   
${\cal F} \equiv {\cal F}(X, d) :=$ the set of $d$-closed subsets of $X$
\eenum 
\edefn 

Note that $B_H(A, r) \subset B_{d_H}(A, r)$. 

\bfactb The function $d_H \in F({\cal S}(X)^2, \widetilde{\IR}_{\geq 0})$ has the following properties. Let $A, A', A'' \in {\cal S}(X)$. 
\benum 
\item (i) \ $A' \subset A''$ \LRA $d_H(A,A') \leq d_H(A,A'')$ \hsh (ii) \ $d_H(x,y) = d(x,y)$ \ \ $(\upfa x,y \in X)$ 

\item[(3)] For $r,s \in \IR_{\geq 0}$ the following holds.
\vskip 0.5mm 
\bit 
\itemI $d_H(A,A') \leq r$ \LLRA $A' \subset B(A, s)$ \ \ $(\upfa s \in \IR_{> r})$ 
\vskip 0.5mm 
\itemII $B(A, r) \subset B_H(A,r) \subset B(A, s)$ \ \ $(\upfa s \in \IR_{> r})$ \hsp $B_H(A, \infty) = B(A, \infty)$ 
\vskip 0.5mm 
\itemiii $d_H(A,A') \leq r$ \LLRA $d_H(A,a') \leq r$ \ \ $(\upfa a' \in A')$ \LLRA $A' \subset B_H(A, r)$ 
\vskip 0.5mm 
\itemiv $d_H(A, B_H(A, r)) \leq r$
\hsp (v) \ $B_H(B_H(A, r), s) \subset B_H(A, r+s)$
\eit 
\vskip 1mm 
\item $d_H(A, A'') \leq d_H(A, A') + d_H(A', A'')$ 
\vskip 0.5mm 
\item 
\bit 
\itemI $d_H(A, A') = 0$ \LLRA $A' \subset B_H(A, 0)$ \LLRA $B_H(A',0) \subset B_H(A, 0)$
\vskip 0.5mm 
\itemII $d_H(A, A) = 0$  
\hsh (iii) \ $d_H(A, A') = 0 = d_H(A', A)$ \LLRA $B_H(A,0) = B_H(A', 0)$ 
\eit 
\vskip 0.5mm 
\item 
\bit 
\itemI $d_H \in {\cal M}_{as}({\cal S}(X))$ \ and \ $\widehat{d_H}$ is a pseudo-metric. 
\vskip 1mm 
\itemII $d_H|_{\cal F} \in {\cal M}_{as}({\cal F})$ is nondegenerate, 
so that $\widehat{d_H}|_{\cal F} = \widehat{d_H|_{\cal F}} \in {\cal M}({\cal F})$, 
where ${\cal F} \equiv {\cal F}(X, d)$. 
\eit 
\eenum 
\efact 

\bfactb Let $d, d' \in {\cal M}_{as}(X)$. The symbols $B, B_H$ and $B', B_H'$ are determined by $d$ and $d'$ respectively. 
\benum 
\item If \ $d' \leq \lambda d + C$ \ for $\lambda \in \IR_{\geq 1}$ and $C \in \IR_{\geq 0}$, then 
\vskip 1mm 
(i) \ $B_H(A,r) \subset B_H'(A, \lambda r + C)$ \ \ \ $(\upfa r \in \widetilde{\IR}_{\geq 0})$ \hsp 
(ii) \ $d_H' \leq \lambda d_H + C$
\vskip 1mm 
\item (i) \ $d \leq_q d'$ \LRA $d_H \leq_q d_H'$ \hsh (ii) \ $d \sim_q d'$ \LRA $d_H \sim_q d_H'$ 
\LRA $\widehat{d}_H \sim_q \widehat{d'}_H$
\eenum
\efact

In the case where $d \in {\cal M}_{as}(X)$ takes discrete values, we do not need to distinguish between $B(A, r)$ and $B_H(A, r)$. 

\bfactb Suppose $J$ is a discrete subset of $\IR_{\geq 0}$ and $d \in {\cal M}_{as}(X, \widetilde{J})$. 
Then, the following holds. 

\benum 
\item 
\bit 
\itemI $0 \in J $ (if $X \neq \emptyset$), since $d(x,x) = 0$ for any $x \in X$.  
\itemII 
For any $r \in \IR_{\geq 0}$ there exists $\lfloor r \rfloor_J := \max J_{\leq r}$, 
since $J_{\leq r} = \{ k \in J \mid k \leq r \}$ is a nonempty finite set of $J$. 
\eit 
\item[] \hspace*{-10mm} Below, let $r \in \IR_{\geq 0}$ and put $n := \lfloor r \rfloor_J $. 
\vskip 1mm 
\item (i) \ $B(A, r) = B(A, n)$ \hsp (ii) \ $d_H(A, A') \leq r$ \LLRA $d_H(A, A') \leq n$ \hsp $(A' \subset X)$
 \vskip 3mm 
 
\item $d_H \in {\cal M}_{as}({\cal S}(X), \widetilde{J})$ : \  
$d_H(A, A') := \left\{ \hspace*{-1mm} \bary[c]{ll}
\min \{ k \in J \mid A' \subset B(A, k) \} & (\upexi  r \in \IR_{\geq 0} \ \ \mbox{s.t.} \ A' \subset B(A, r)) \\[2mm]
\ \infty & (\mbox{otherwise}) 
\eary \right.$
\vskip 2mm 
\item 
(i) \ $B_H(A,r) = B_H(A, n)$ 
\hsh (ii) \ $B_H(A,r) = B(A, r)$
\eenum 
\efact 

\section{Metric on the power set of a group induced from the word length}  

\subsection{Basics on operations for subsets of a group} \mbox{}

Suppose $G$ is a group and $e \in G$ is the identity element. Let $G^\times := G - \{ e \}$. 

\bdefnb\label{def_norm}
\benum 
\item A norm on $G$ means a function $\nu : G \to \widetilde{\IR}_{\geq 0}$ which satisfies the following conditions : 
\bit 
\itemI $\nu(x) = 0$ \LLRA $x = e$ 
\itemII $\nu(x^{-1}) = \nu(x)$ \ $(\upfa x \in G)$
\itemiii $\nu(xy) \leq \nu(x) + \nu(y)$ \ $(\upfa x,y \in G)$ 
\eit 
\item A function $\nu : G \to \widetilde{\IR}_{\geq 0}$ is conjugation-invariant if it satisfies the following condition : 
\bit 
\itemiv $\nu(axa^{-1}) = \nu(x)$ \ $(\upfa x, a \in G)$ 
\eit 
\eenum 
\edefn

\bnotb For $S, T \subset G$ we use the following notations.
\benum
\item $ST := \{ xy \mid x \in S, y \in T \} \subset G$
\item $S^n := \{ x_1 \cdots x_n \mid x_i \in S \ (i \in [n]) \}$ \ \ $(n \in \IZ_{\geq 1})$ 
\hsh \hsf $S^0 := \{ e \}$ \hsh \hsf $S^{-1} = \{ x^{-1} \mid x \in S \}$ 
\item $S^{\leq n} = \bigcup_{k = 0}^n S^k$ \ \ $(n \in \IZ_{\geq 1})$ \hsp $S^\infty = \bigcup_{n = 0}^\infty S^n$
\item (the conjugations)
\bit 
\itemI $x^a := a^{-1}xa$, \ ${}^ax := axa^{-1}$ \ $(x, a \in G)$ \hsp (ii) \ 
$S^A := \{ x^a \mid x \in S, a \in A \}$ \ $(A \subset G)$ 
\eit 
\eenum 
\enot 

\bnotb For $S \subset G$ we use the following notations.
\benum 
\item 
\bit 
\itemI $\langle S \rangle$ : the subgroup of $G$ generated by $S$ 
\itemII $\langle\langle S \rangle\rangle$ : the normal subgroup of $G$ normally generated by $S$ 
\eit 
\item $S^\pm := S \cup S^{-1}$ \ (the symmetrization of $S$)
\item 
\bit 
\itemI $C(S) := S^G$ \ (the conjugation closure of $S$) 
\itemII $C_S := (S^\pm)^G = (S^G)^{\pm}$ \ (the symmetrized conjugation closure of $S$)
\eit 
\vskip 1mm 
\item 
\bit 
\itemI $S$ is symmetric. \LLRAdefn $S = S^{-1}$
\itemII $S$ is conjugation-invariant. \LLRAdefn $S^a \subset S$ \ $(\upfa a \in G)$ \LLRAequiv $S^a = S$ \ $(\upfa a \in G)$
\eit 
\eenum 
\enot

\bnot\label{not_S_p} For subsets of ${\cal S}(G)$ we use the following notations. 
\bit 
\item[{[I]}] For any condition $P$ on subsets $S$ of $G$, we set ${\cal S}_{P}(G) := \{ S \in {\cal S}(G) \mid S : P \}$.  \\
Examples of a condition $P$ : 
\benum 
\item 
\bit 
\itemI $f$ : \ $|S| < \infty$ \ (a finite set) \hsp (ii) \ $g$ : \ $S^\infty = G$ \hsp (iii) \  $s$ : \ $S^{-1} = S$ \ (symmetric) 
\eit 
\item
\bit 
\itemI $c$ : \ $S^G = S$ \ (conjugation-invariant) \LLRA $S = A^G$ \ $(\upexi A \in {\cal S}(G))$ 
\itemII $fc$ : \ $S = A^G$ \ $(\upexi A \in {\cal S}_f(G))$ 
\hsp (iii) \ $s, c := s \,+\, c$ \LLRA $S = C_A$ \ $(\upexi A \in {\cal S}(G))$  
\eit 
\item 
\bit 
\itemI $fg := g \,+\, s \,+\, f$ 
\hsp (ii) \ $ng := g \,+\, s \,+\, c$ 
\itemiii $nfg := g \,+\, s \,+\, fc$ \LLRA $S^\infty = G$, \ $S = C_A$ \ $(\upexi A \in {\cal S}_f(G))$ 
\eit 
\eenum 
\vskip 1mm 
\item[{[II]}] To avoid the degeneracy of induced metrics on ${\cal S}(G)$ (cf. Fact~\ref{fact_nu_H}\,(1)), 
we need some care on the identity element $e$.
For this reason we introduce the following notations. 
\benum 
\item For ${\cal T} \subset {\cal S}(G)$ we set ${\cal T}^\ast := {\cal T} - \{ \emptyset, \{ e \} \}$.  

\item ${\cal S}(G) = {\cal S}'(G) \coprod {\cal S}''(G)$ : \hsh  
${\cal S}'(G) := \{ S \in {\cal S}(G) \mid e \not\in S \}$, \hsh ${\cal S}''(G) := \{ S \in {\cal S}(G) \mid e \in S \}$ 
\bit 
\itemI ${\cal S}'_P(G) := {\cal S}'(G) \cap {\cal S}_P(G)$, \hsh ${\cal S}''_P(G) := {\cal S}''(G) \cap {\cal S}_P(G)$
\eit  

\item $\widetilde{\cal S}(G) := {\cal S}(G)/\sim$ : 
\bit 
\itemI $\sim$ is the equivalence relation on ${\cal S}(G)$ generated by $S \sim S \cup \{ e \}$ \ $(S \in {\cal S}(G))$. 
\itemII $S \sim T$ \LLRA $S - \{ e \} = T - \{ e \}$ \hsh $(S, T \in {\cal S}(G))$ 
\itemiii The equivalence class of $S \in {\cal S}(G)$ is denoted by $[S]$. 
\itemiv $\widetilde{\cal S}_P(G) := \{ [S] \in \widetilde{\cal S}(G) \mid S \in {\cal S}_P(G) \}$ 
\eit 
\item There exist canonical bijections \ \ 
$\bary[t]{@{}c@{ \ }c@{ \ }c@{ \ }c@{ \ }cl@{}}
{\cal S}'(G) & \approx & {\cal S}''(G) & \approx & \widetilde{\cal S}(G) & \\[1mm] 
S & & S \cup \{ e \} & & [S] & \hspace*{-4mm} = \{ S, S \cup \{ e \} \}.
\eary$
\eenum 
\eit 
\enot

Below we prefer the subset ${\cal S}''(G)$ since related expressions take a simple form (cf. Fact~\ref{fact_formula}\,(4)(iv)). 
Next we list up basic properties of subsets of $G$ under various operations. 

\bfact\label{fact_formula} Let $G \supset S, S_1, S_2, T, T_1, T_2, U, S_\lambda, T_\mu$ \ $(\lambda \in \Lambda, \mu \in M)$ and $G \ni a,b$.  
\benum 
\item  
\bit 
\itemI $S_1 \subset S_2$, \ $T_1 \subset T_2$ \LRA $S_1T_1 \subset S_2T_2$, \ \ \ $S_1^{T_1} \subset S_2^{T_2}$
\hsp (ii) \ $eS = S = Se$
\vskip 0.5mm 
\itemiii $(ST)U = S(TU)$ \hsh (iv) \ $(S^T)^U = S^{TU}$ 
\itemv $(S \cup T)^{-1} = S^{-1} \cup T^{-1}$ \hsh $(ST)^{-1} = T^{-1}S^{-1}$ \hsh $(S^T)^{-1} = (S^{-1})^T$ 
\item[(vi)] $TU \cup TU^{-1} \subset T$ \LRA $S^TU = US^T$ 
\bit 
\item[\bec] $S^TU \ni x^ay = yy^{-1}x^ay = y (x^a)^y = y x^{ay} \in US^T$ \\
$US^T \ni y x^a = yx^ay^{-1}y = (x^a)^{y^{-1}}y = x^{ay^{-1}}y \in S^TU$ 
\eit  
\eit 
\item $(S \cup T)^n = \bigcup_{I + J = [n]} S^I T^J$ \hsh ($S^I T^J = X_1 \cdots X_n$ : \ $X_k = S$ if $k \in I$, $X_k = T$ if $k \in J$) 
\bit 
\itemI $ST = TS$ \LRA $(S \cup T)^n = \bigcup \{ S^k T^\ell \mid k, \ell \in \IZ_{\geq 0}, k + \ell = n \}$
\hsh $(ST)^n = S^nT^n$ 
\itemII $(S \cup T)^n = S^n \cup S^{\leq n-1}T$ \ \ if $T \vartriangleleft G$. \hspp $(S^{\leq -1} := \emptyset)$ 
\eit 

\item 
\bit 
\itemI $S = \bigcup_{\lambda \in \Lambda} S_\lambda$, $T = \bigcup_{\mu \in M} T_\mu$ \\[1mm] 
\LRA $ST = \bigcup_{(\lambda, \mu) \in \Lambda \times M} S_\lambda T_\mu$, \ 
$S^T = \bigcup_{(\lambda, \mu) \in \Lambda \times M} S_\lambda^{T_\mu}$, \ $S^{-1} = \bigcup_{\lambda \in \Lambda} S_\lambda^{-1}$
\vskip 1mm 
\itemII $S = \bigcap_{\lambda \in \Lambda} S_\lambda$ \LRA $S^{-1} = \bigcap_{\lambda \in \Lambda} S_\lambda^{-1}$, \ 
$S^T \subset \bigcap_{\lambda \in \Lambda} S_\lambda^T$
\eit 
\vskip 1mm 
\item 
\bit 
\itemI $S^k S^\ell = S^{k+\ell}$ \hsh $(S^k)^\ell = S^{k \ell}$ \hsh $(k,\ell \in \IZ_{\geq 0})$ \hsp $S^\infty S^\infty = S^\infty$
\itemII $(S^n)^{-1} = (S^{-1})^n$ \hsh $(S^{\leq n})^{-1} = (S^{-1})^{\leq n}$ \hsh $(n \in \widetilde{\IZ}_{\geq 0})$ 
\vskip 0.5mm 
\itemiii $S^{\leq m} S^{\leq n} = S^{\leq m+n}$ \hsh $(S ^{\leq m})^{\leq n} \subset S^{\leq mn}$ \hsh $(m, n \in \IZ_{\geq 0})$  
\itemiv $e \in S$ \LRA 
$S^n = S^{\leq n}$ \hsh $(n \in \IZ_{\geq 0})$ \hsp (v) \ $(S \cup \{ e \})^n = (S \cup \{ e \})^{\leq n} = S^{\leq n}$ \hsh $(n \in \IZ_{\geq 0})$
\vskip 0.5mm 
\item[(vi)] $(x^a)^b = x^{ab}$ \hsh $(xy)^a = x^ay^a$ \hsh $(x^{-1})^a = (x^a)^{-1}$ \hsh $(S \cup T)^a = S^a \cup T^a$ 
\hsh $(S \cap T)^a = S^a \cap T^a$ \\ 
$(ST)^a = S^aT^a$ \hsh $(S^a)^b = S^{ab}$ \hsh $(S^{-1})^a = (S^a)^{-1}$ \hsh  \\
 $(S^n)^a = (S^a)^n$ \hsh 
$(S^{\leq n})^a = (S^a)^{\leq n}$ \hsh $(S^{\leq n})^T \subset (S^T)^{\leq n}$ \hsh  $(n \in \widetilde{\IZ}_{\geq 0})$
\eit

\item 
\bit 
\itemI $S^\pm$ is symmetric. \hsh $(S \cup T)^{\pm} = S^{\pm} \cup T^{\pm}$ 
\itemII $S, T$ are symmetric.  \LRA $S^{-1}$, $S \cup T$, $S \cap T$ are symmetric. \hsh $(ST)^{-1} = TS$ 
\itemiii $S_\lambda$ is symmetric \ $(\upfa \lambda \in \Lambda)$. \LRA 
$\bigcap_{\lambda \in \Lambda} S_\lambda$ is symmetric. 
\eit 

\item 
\bit 
\itemI ${}^T\! S = S^{T^{-1}}$ \hsh (ii) \ $ST \subset TS^T$, \ \ $TS \subset {}^T\!S\, T = S^{T^{-1}}T$ 
\itemiii $S$ is conjugation-invariant. \LRA ${}^T\!S = S^T = S$, \ \ $ST = TS$ 
\itemiv $S, T$ are conjugation-invariant. \LRA $S^{-1}$, \ $S \cup T$, \ $S \cap T$, \ $ST$ are conjugation-invariant.  
\itemv $S_\lambda$ is conjugation-invariant \ $(\upfa \lambda \in \Lambda)$. \LRA 
$S = \bigcap_{\lambda \in \Lambda} S_\lambda$ is conjugation-invariant. 
\item[(vi)] $S, T$ are symmetric, conjugation-invariant. \\
\hsp \LRA $S^{-1}$, $S \cup T$, \ $ST$ are symmetric, conjugation-invariant.  
\eit 
\item
\bit 
\itemI $C(S)$ is conjugation-invariant. \hsh (ii) \ $C_S$ is symmetric and conjugation-invariant.  
\itemiii $C_{S \cup T} = C_S \cup C_T$ \hsh $C_{S \cup \{ e \}} = C_S \cup \{ e \}$ \hsp 
(iv) \ $S \subset G^\times$ \LRA $C_S \subset G^\times$ 
\eit 
\vskip 1mm 
\item 
\bit 
\itemI $S^\infty$ is a submonoid of $G$ generated by $S$. 
\itemII $\langle S \rangle = \langle S^\pm \rangle = (S^\pm)^\infty$ \hsh 
(iii) \ $\langle\langle S \rangle\rangle = \langle C(S) \rangle = (C_S)^\infty$ 
\hsp $\circ$ \ $\emptyset^0 = \emptyset^\infty = \langle \emptyset \rangle = \langle\langle \emptyset  \rangle\rangle = \{ e \}$ 
\eit 
\eenum 
\efact 

\subsection{Metrics on the power set of a group} \mbox{} 

\bdefn\label{defn_d_S} Let $S \in {\cal S}(G)$.  
\benum 
\item the word metric on $G$ with respect to $S$ : \\[1mm] 
\hsp $d_S : G \times G \lra \widetilde{\IZ}_{\geq 0}$ : \ 
$d_S(x, y) := 
\left\{ \hspace*{-1mm} \bary[c]{ll}
\min \{ n \in \IZ_{\geq 0} \mid y \in x S^n \} & (y \in xS^\infty) \\[1mm]
\infty & (x \in G - xS^\infty) 
\eary \right.$ 
\vskip 2mm 

\item the Hausdorff metric induced from $d_S$ : \\[1mm] 
\hsp $(d_S)_H : {\cal S}(G) \times {\cal S}(G) \lra \widetilde{\IZ}_{\geq 0}$ : \\[2mm]  
\hspp    
$(d_S)_H(A, A') := \left\{ \hspace*{-1mm} \bary[c]{ll}
\min \{ n \in \IZ_{\geq 0} \mid A' \subset AS^{\leq n} \} & (\upexi  n \in \IZ_{\geq 0} \ \ \mbox{s.t.} \ A' \subset AS^{\leq n}) \\[2mm]
\ \infty & (\mbox{otherwise}) 
\eary \right.$
\vskip 2mm 
\bit 
\itemI For $n \in \IZ_{\geq 0}$ \hsh $(d_S)_H(A, A') \leq n$ \LLRA $A' \subset A S^{\leq n}$ 
\eit 
\eenum 
\edefn
\vskip 1mm 

\bdefnb Let $S \in {\cal S}(G)$.
\benum 
\item the word length norm on $G$ with respect to $S$ : \\[1mm] 
\hsh $\nu_S := d_S(e, \ast) : G \lra \widetilde{\IZ}_{\geq 0}$ : \ 
$\nu_S(x) = d_S(e, x) = 
\left\{ \hspace*{-1mm} \bary[c]{ll}
\min \{ n \in \IZ_{\geq 0} \mid x \in S^n \} & (x \in S^\infty) \\[1mm]
\infty & (x \in G - S^\infty)  
\eary \right.$ 
\vskip 2mm 
\bit 
\item[]  
\bit 
\itemI For $n \in \IZ_{\geq 0}$ \hsh $\nu_S(x) \leq n$ \LLRA $x \in S^{\leq n}$
\eit 
\eit 
\vskip 2mm 

\item the Hausdorff-like extension of $\nu_S$ : \\[1mm]  
\hsh $(\nu_S)_H := (d_S)_H(e, \ast) : {\cal S}(G) \lra \widetilde{\IZ}_{\geq 0}$ : \ $(\nu_S)_H(A) = (d_S)_H(e, A) = \sup\,\nu_S(A)$ \\[2mm] 
\hspace*{90mm} ($\sup$ is taken in $\widetilde{\IR}_{\geq 0}$ so that $\sup \emptyset = 0$.)  \\[1mm] 
\hsh $\rho_S :=\log\, (\nu_S)_H : {\cal S}(G)^\ast \lra \log\,(\widetilde{\IZ}_{\geq 1})\ \subset \ \widetilde{\IR}_{\geq 0}$ 
\bit 
\item[]
\bit 
\itemI For $n \in \IZ_{\geq 0}$ \hsh $(\nu_S)_H(A) \leq n$ \LLRA $A \subset S^{\leq n}$ 
\itemII $(\nu_S)_H(A) = 0$ \LLRA $A \subset \{ e \}$ \LLRA $A = \{ e \}$ or $\emptyset$ 
\eit 
\eit 
\vskip 1mm 

\item[] $\circ$ \ If there is no ambiguity we omit the subscript $H$ and write $(\nu_S)_H$ simply as $\nu_S$. 
\eenum 
\edefn

\bdefnb 
\benum 
\item $\nu_H  : {\cal S}(G) \times {\cal S}(G) \lra \widetilde{\IZ}_{\geq 0}$ : \ $\nu_H(S, T) := (\nu_S)_H(T)$ \\[1mm]
$\rho : = \log \nu_H : {\cal S}(G)^\ast \times {\cal S}(G)^\ast \lra \widetilde{J} \ \subset \ \widetilde{\IR}_{\geq 0}$ : \ $\rho(S,T) = \rho_S(T)$ 
\vskip 1mm 
\item $\widehat{\rho} := \log \widehat{\nu_H} : {\cal S}(G)^\ast \times {\cal S}(G)^\ast 
\lra \widetilde{J} \ \subset \ \widetilde{\IR}_{\geq 0}$ \\[1mm]  
\hsh $\widehat{\rho}(S,T) = \max\,\{ \rho_S(T), \rho_T(S) \} = \log \max \{ (\nu_S)_H(T), (\nu_T)_H(S) \}
= \log\,\sup\,(\nu_S(T) \cup \nu_T(S))$ 
\eenum 
\edefn  

\bfact\label{fact_d-nu} Let $S \in {\cal S}(G)$. 
\benum 
\item $d_S \in {\cal M}_{as}(G, \widetilde{\IZ}_{\geq 0})$
\vskip 1mm 
\bit 
\itemI $d_S(x,y) = 0$ \LLRA $x = y$ \hsh (ii) \ $d_S(x,y) = d_{S^{-1}}(y,x)$ \ \ $(\upfa x,y \in G)$ 
\vskip 0.5mm 
\itemiii $d_S(ax, ay) = d_S(x,y)$ \ \ $(\upfa a, x,y \in G)$ 
\hsh (iv) \ $d_{S^a}(x^a, y^a) = d_S(x,y)$ \ \ $(\upfa a, x,y \in G)$ 
\eit 
\vskip 1mm 
\item $\nu_S$ satisfies the conditions (i), (iii) in Definition~\ref{def_norm} 
\vskip 1mm 
\bit 
\itemI $d_S(x, y) = \nu_S(x^{-1}y)$ \ \ $(\upfa x,y \in G)$ \hsh (ii) \  $\nu_{S^a}(x^a) = \nu_S(x)$ \ \ $(\upfa a, x \in G)$ 
\vskip 1mm 
\itemiii $\nu_S^{-1}(0) = \{ e \}$, \ $\nu_S^{-1}(1) = S - \{ e \}$
\eit 
\vskip 1mm 
\item If $S = S^{-1}$, then \ (i) \ $d_S \in {\cal M}(G, \widetilde{\IZ}_{\geq 0})$ \ \ \ (ii) \ $\nu_S$ is a norm on $G$. 
\item If $S$ is conjugation invariant, then $d_S$ and $\nu_S$ are conjugate-invariant. 
\eenum 
\efact 

\bfactb\label{fact_nu_H}  
\benum 
\item Let $S, T, A \subset G$. 
\bit 
\itemI $(\nu_S)_H(A) = 1$ \LLRA $A \not\subset \{ e \}$ and $A \subset S \cup \{ e \}$ 
\itemII $(\nu_S)_H(T) = (\nu_T)_H(S) = 1$ \LLRA $S - \{ e \} = T - \{ e \} \neq \emptyset$
\itemiii $(\nu_S)_H(A \cup \{ e \}) = (\nu_S)_H(A)$, \ \ $(\nu_{S \cup \{ e \}})_H(A) = (\nu_S)_H(A)$ 
\eit 
\item $\nu_H \in {\cal M}_{m, as}({\cal S}(G)^\ast)$ and $\rho \in {\cal M}_{as}({\cal S}(G)^\ast)$. 

\item $\widehat{\rho} \in {\cal M}({\cal S}''(G)^\ast)$ 
\item $\nu_H$ induces the function : \ 
$\widetilde{\nu}_H : \widetilde{\cal S}(G) \times \widetilde{\cal S}(G) \lra \widetilde{Z}_{\geq 0}$ : $\widetilde{\nu}_H([S], [T]) = \nu_H(S,T)$. 
\vskip 1mm 
\item[] $\widetilde{\rho} := \log \,\widetilde{\nu}_H : \widetilde{\cal S}(G) \times \widetilde{\cal S}(G) \lra \log \,\widetilde{Z}_{\geq 0}$
\vskip 1mm 
\bit 
\itemI We have natural isometric bijections \ \ 
$\bary[t]{@{}c@{ \ }c@{ \ }c@{ \ }c@{ \ }c@{}}
({\cal S}'(G), \nu_H) & \approx & ({\cal S}''(G), \nu_H) & \approx & (\widetilde{\cal S}(G), \widetilde{\nu}_H) \\[1mm] 
S & & S \cup \{ e \} & & [S] 
\eary.$ 
\vskip 1mm 
\itemII $\nu_H$ is nondegenerate on ${\cal S}'(G)^\ast$ and ${\cal S}''(G)^\ast$ and 
$\widetilde{\nu}_H$ is nondegenerate on $\widetilde{\cal S}(G)^\ast$.
\eit 
\item For $S, T \in {\cal S}''(G)$ \hsh (i) \ $S \cup T \subset ST \subset (S \cup T)^2$ \hsh 
(ii) \ $\nu_H(S\cup T, ST) \leq 2$, \ $\nu_H(ST, S \cup T) \leq 1$
\eenum 
\efact 

\vskip 1mm 
\bfactb There exists an isometric bijection \\[1mm] 
\hspace*{30mm} $\phi : (F(G - A, J), \mu) \lra (F_L(G^2 - \widetilde{A}, J), \mu)$ : $f \lmt \widetilde{f} := f \eta|$. \\[1mm] 
\hsp Here, \btab[t]{cl} 
(i) & $\eta : G^2 \, \to\hspace{-4.5mm}\to \, G$ : $\eta(x,y) = x^{-1}y$, \ \ \ $J \subset \widetilde{\IR}_{\geq 0}$ \\[2mm] 
(ii) & $A \subset G$, \ \ \ $\widetilde{A} := \eta^{-1}(A)$, \ \ \ $\eta| : G^2 - \widetilde{A} \,\lra\hspace{-5.5mm}\lra \, G - A$ \\[2mm] 
(iii) & $F_L(G^2 - \widetilde{A}, J) := \{ g \in F(G^2 - \widetilde{A}, J) \mid g(ax, ay) = g(x,y)$ \ $(\upfa a \in G, (x,y) \in G^2 - \widetilde{A}\,) \}$ \\[2mm] 
(iv) & The metric $\mu$ is defined in Notation~\ref{not_F}\,(1). \\[2mm] 
(v) & $\widetilde{f}(x,y) = f(x^{-1}y)$, \ $f(x) = \widetilde{f}(e,x)$ \ $(x,y \in G)$ 
\etab
\efact 
\vskip 2mm 

\bpropb There exists a natural isometric injection \\[0mm] 
\hspp \hsh $\bary[t]{@{}c@{ \ }c@{ \ }c@{ \ }c@{ \ }c@{ \ }c@{}} 
& & & & \phi & \\[-1mm] 
\zeta : \, & ({\cal S}'(G)^\ast, \widehat{\rho}) & 
\xymatrix@M+1pt{ 
\ar@{^{(}->}[r] &  
} & (F(G^\times, \widetilde{\IZ}_{\geq 0}), \widehat{\mu}) & \ \approx \ & (F_L(G^2 - \Delta(G), \widetilde{\IZ}_{\geq 0}), \widehat{\mu}) \\[0.5mm]
& S & & \log \nu_S & & \log d_S 
\eary$ 
\bit 
\itemC $\zeta$ is anti-isometric with respect to $\rho$ and $\mu$ \ 
(i.e., $\rho(S,T) = \mu(\log \nu_T, \log \nu_S)$ \ $(\upfa S, T \in {\cal S}'(G)^\ast$). 
\eit 
\eprop 
\vspace*{-4mm} 
\bpfb 
\benum 
\item Consider the map $\psi : F(G, \widetilde{\IZ}_{\geq 0}) \lra {\cal S}'(G)$ : $\psi(f) = f^{-1}(0)$. 
Since $\psi \zeta = \id$ by Fact~\ref{fact_d-nu}\,(2)(iii), 
it follows that $\zeta$ is injective. 
\vskip 1mm 
\item $\nu_H(S,T) \equiv (\nu_S)_H(T) = \lambda(\nu_T, \nu_S)$ \ \ $(\upfa S, T \in {\cal S}(G))$. \\[1mm] 
Hence, the map \ 
$\bary[t]{@{}c@{ \ }c@{ \ }c@{}} 
({\cal S}'(G), \nu_H)  & \lra & (F(G^\times, \widetilde{\IZ}_{\geq 1}), \lambda) \\[0mm] 
S & & \nu_S 
\eary$ \ is an anti-isometric injection. \\[1mm] 
This follows from the following observation : 
\bit 
\item[] For $n \in \IZ_{\geq 1}$ \hsh 
$(\nu_S)_H(T) \leq n$ \LLRA $T \subset S^{\leq n}$ \ $\relsb{}{\Llra}{\mbox{$(\ast)$}}$ \ $\nu_S \leq n\, \nu_T$ 
\LLRA $\lambda(\nu_T, \nu_S) \leq n$
\eit 
\vskip 1mm 
\bit 
\item[$(\ast)$]
$(\Lra)$ \ For  any $x \in G$ we have to show that \ $\nu_S(x) \leq n \,\nu_T(x)$. Let $m := \nu_T(x)$. \\
\hsp \hsf If $m = \infty$, then the assertion is trivial. \\
\hsp \hsf If $m < \infty$, then $x \in T^m \subset (S^{\leq n})^m \subset S^{\leq mn}$.
This implies that $\nu_S(x) \leq mn = n\,\nu_T(x)$. 
\item[] $(\Lla)$ \ For any $x \in T$ it follows that $\nu_T(x) \leq 1$ and $\nu_S(x) \leq n \,\nu_T(x) \leq n$, so that $x \in S^{\leq n}$.   
\eit 
\vskip 1mm 
\item For any $S, T \in {\cal S}(G)^\ast$ it follows that \\
\hsp $\rho(S,T) \equiv \rho_S(T) = \log\, (\nu_S)_H(T) = \log \lambda(\nu_T, \nu_S) = \mu(\log \nu_T, \log \nu_S)$ \ \ and \\[1mm] 
\hsp $\widehat{\rho}(S,T) = \widehat{\mu}(\log \nu_S, \log \nu_T)$. 
\eenum 
\vspace*{-7mm} 
\epf

Here we have some examples of subspaces ${\cal X} \equiv ({\cal X}, \widehat{\rho}\,|_{{\cal X}})$ of 
$({\cal S}'(G), \widehat{\rho}\,)$ with $\widehat{\rho}\,|_{{\cal X}} < \infty$

\bexp\label{exp_finite} \mbox{} 
\benum 
\item Suppose $G$ is finitely generated and $n := {\rm rank}\,G \geq 1$. \\ 
\hsh ${\cal S}_{fg}'(G) := \{ S \in {\cal S}'(G) \mid S = S^{-1}, \ G = S^\infty, \ |S| < \infty \}$ \\[1mm] 
\hsh ${\cal S}_{fg}'(G)_n := \{ S \in {\cal S}'_{fg}(G) \mid |S| = n \}$ 
\vskip 1mm 
\item Suppose $G$ is finitely normally generated and $n := \mbox{n-rank}\,G \geq 1$. \\ 
\hsh ${\cal S}_{nfg}'(G)^\ast := \{ S \in {\cal S}'(G)^\ast \mid G = S^\infty, \ S = C_A \mbox{ for some $A \subset G^\times$ \ with \ $|A| < \infty$} \}$ \\[1mm]
\hsh ${\cal S}_{nfg}'(G)^\ast_n := \{ S \in {\cal S}'(G)^\ast \mid G = S^\infty, \ S = C_A \mbox{ for some $A \subset G^\times$ \ with \ $|A| = n$} \}$  
\vskip 1mm 

\item[] $\circ$ \ $(\nu_S)_H(T) \in \IR_{\geq 1}$ \ $(\upfa S, T \in {\cal S}'_{nfg}(G)^\ast)$ \\
\ (Proof) 
\bit 
\itemi Since $S = C_A$ for some  $A \subset G^\times$, it follows that 
$S$ is symmetric and conjutgation invariant. 
Hence, $\nu_S$ is a conjutgation invariant norm on $G$ (cf. Fact~\ref{fact_d-nu}\,(3), (4)).   

\itemii Since $G = S^\infty$ and $T = C_B$ for some  $B \subset G^\times$ \ with \ $|B| < \infty$, it follows tnat \\
\hsf \ $\nu_S(T) = \nu_S((B^\pm)^G) = \nu_S(B^\pm) = \nu_S(B)$ \ and \ 
$(\nu_S)_H(T) = \sup\,\nu_S(T) = \max\, \nu_S(B) < \infty$.  
\eit 
\eenum 
\eexp 

\section{The maps on the power sets induced from group homomorphisms} 

Suppose $f : G \to H$ is a group homomorphism between groups. Let $K := {\rm Ker}\,f$. 
We are interested in the following general problem. 

\begin{problem} Describe the relation among the metric spaces $({\cal S}''(G)^\ast, \nu_H^G)$, $({\cal S}''(H)^\ast, \nu_H^H)$
and $({\cal S}''(K)^\ast, \nu_H^K)$. 
\end{problem}

\subsection{Generalities} \mbox{} 

First we list some general properties. 

\bfactb\label{fact_f+f^{-1}}
\benum 
\item[] \hspace*{-10mm} {[I]} \ Let $G \supset S, A, B, A_1, \cdots,  A_n$. 
\item 
\bit 
\itemI $f(AB) = f(A)f(B)$ 
\hsh \ (ii) \ $f(A_1 \cdots A_n) = f(A_1) \cdots f(A_n)$ 
\hsh \ (iii) \ $f(A^B) = f(A)^{f(B)}$ 
\itemiv $f(A^{-1}) = f(A)^{-1}$ 
\eit 
\vskip 1mm 
\item 
\bit 
\itemI $f(S^n) = f(S)^n$, \ $f(S^{\leq n}) = f(S)^{\leq n}$ \ \ $(\upfa n \in \IZ_{\geq 0})$ \hsh (ii) \ $f(S^\infty) = f(S)^\infty$  
\eit 
\vskip 1mm 

\item Suppose $f$ is surjective. 
\bit 
\itemI $f(C_A) = C_{f(A)}$
\itemII The following implication holds for 
the conditions $P = f, \ g,\ s, \ c, \ fc, \ fg, \ ng, \ nfg$ : \\
\hsh $S \in {\cal S}''_P(G)$ \LRA $f(S) \in {\cal S}''_P(H)$ 
\eit 
\vskip 1mm 
\item[] \hspace*{-10mm} \!{[II]} \ Suppose $f$ is surjective. Let $H \supset S, A, B, A_1, \cdots,  A_n$
\vskip 0.5mm 
\item[(1)] 
\bit 
\itemI  $f^{-1}(AB) = f^{-1}(A)f^{-1}(B)$ \hsp (ii) \ $f^{-1}(A_1\cdots A_n) = f^{-1}(A_1)\cdots f^{-1}(A_n)$ \\[1mm] 
\ \ $\circ$ \ $f^{-1}(A) K = f^{-1}(A)$ 
\eit 
\vskip 1mm 
\item[(2)]  
\bit 
\itemI $f^{-1}(S^n) = (f^{-1}(S))^n$ \ \ $(\upfa n \in \IZ_{\geq 1})$ 
\vskip 2mm 
\itemII $f^{-1}(S^{\leq n}) = (f^{-1}(S))^{\leq n} \cup K
= \left\{ \hspace*{-1mm} \bary[c]{ll}
K & (n=0) \\[2mm]
(f^{-1}(S) \cup K)^{\leq n} & (n \in \IZ_{\geq 1}) 
\eary \right.$ 
\vskip 2mm  

\itemiii $f^{-1}(A S^{\leq n}) = f^{-1}(A)(f^{-1}(S))^{\leq n}$ \hsh $(\upfa n \in \IZ_{\geq 0})$
\vskip 1mm 
\itemiv $f^{-1}(S^{-1}) = (f^{-1}(S))^{-1}$  
\hsp (v) \ $f^{-1}(S^\infty) = (f^{-1}(S))^\infty \cup K$ 
\eit 
\eenum 
\efact 

\bnotb The group homomorphism $f : G \to H$ induces the following maps. \\[1mm]  
\hsp \hsh 
\btab[t]{ll}
$f_\ast : {\cal S}(G) \lra {\cal S}(H)$ : $f_\ast(S) = f(S)$ &  
$\widetilde{f}_\ast : \widetilde{\cal S}(G) \lra \widetilde{\cal S}(H)$ : $\widetilde{f}_\ast([S]) = [f(S)]$ \\[2mm]
$f^\ast : {\cal S}(H) \lra {\cal S}(G)$ : $f^\ast(S) = f^{-1}(S)$ \hsh & 
\etab 
\vskip 2mm 
Fact~\ref{fact_f+f^{-1}} indicates that these maps are compatible with various operations in ${\cal S}(G)$ and ${\cal S}(H)$. 
\enot 

\bfact\label{fact_f_ast} Suppose $f : G \to H$ is a group homomorphism. 
\benum 
\item Let $S \subset G$.  
\bit 
\itemI $d_{f(S)}(f(x), f(y)) \leq d_S(x,y)$ \hsh $\nu_{f(S)}(f(x)) \leq \nu_S(x)$ \hsh  $(\upfa x,y \in G)$ %$(\upfa x \in G)$
\vskip 1mm 
\itemII $(d_{f(S)})_H(f(A), f(B)) \leq (d_S)_H(A,B)$ \hsh $(\upfa A, B \subset G)$ \\
$(\nu_{f(S)})_H(f(A)) \leq (\nu_S)_H(A)$ \hsh $\nu_H(f(S), f(T)) \leq \nu_H(S, T)$ 
 
\vskip 1mm 
\itemiii $f_\ast : ({\cal S}(G), (d_S)_H) \lra ({\cal S}(H), (d_{f(S)})_H)$ is a 1-Lipschitz map. 
\eit 
\vskip 1mm 
\item 
$f_\ast : ({\cal S}(G)^\ast - {\cal S}(K), \nu_H\,) \lra ({\cal S}(H)^\ast, \nu_H\,)$ is a 1-Lipschitz map. \\
$f_\ast : ({\cal S}(G)^\ast - {\cal S}(K), \widehat{\rho}\,) \lra ({\cal S}(H)^\ast, \widehat{\rho}\,)$ is a 1-Lipschitz map.  
\vskip 1mm 
\item $\widetilde{f}_\ast : (\widetilde{\cal S}(G)^\ast, \widehat{\rho}\,) \lra (\widetilde{\cal S}(H)^\ast, \widehat{\rho}\,)$ is a 1-Lipschitz map. 
\vskip 1mm 
\item[] \hspace*{-10mm} Suppose $f$ is surjective. 
\item Let $S \subset H$. 
\bit 
\itemI $(d_{f^{-1}(S)})_H(f^{-1}(A), f^{-1}(B)) = (d_S)_H(A,B)$ \hsh $(\upfa A, B \subset H)$ \\
$(\nu_{f^{-1}(S)})_H(f^{-1}(A)) = (\nu_S)_H(A)$ \hsh $(\upfa A \in {\cal S}(H)^\ast)$ 
\vskip 1mm 
\itemII $f^\ast : ({\cal S}(H), (d_S)_H) \lra ({\cal S}(G), (d_{f^{-1}(S)})_H)$ is an isometric injection. \hsh $f_\ast f^\ast = \id$  
\eit 
\vskip 0.5mm 
\item $f^\ast : ({\cal S}''(H)^\ast, \nu_H) \lra ({\cal S}''(G)^\ast, \nu_H)$ is an isometric injection. \\
$f^\ast : ({\cal S}''(H)^\ast, \widehat{\rho}\,) \lra ({\cal S}''(G)^\ast, \widehat{\rho}\,)$ is an isometric injection. 
\eenum 
\efact 

\bprop\label{prop_f^ast+f_ast} Suppose $f : G \to H$ is a group epimorphism and $K := {\rm Ker}\,f$. 
\benum 
\item $f_\ast : ({\cal S}''(G)^\ast - {\cal S}(K), \nu_H) \lra ({\cal S}''(H)^\ast, \nu_H)$ is a 1-Lipschitz map. 
\item $f^\ast : ({\cal S}''(H)^\ast, \nu_H) \lra ({\cal S}''(G)^\ast - {\cal S}(K), \nu_H)$ is an isometric section of $f_\ast$. 
\hsh $(f_\ast f^\ast = \id)$  
\item $\widehat{\nu}_H(T, f^\ast f_\ast(T)) \leq \nu_H(T, K) + 1$ \hsh $(\upfa T \in {\cal S}''(G)^\ast - {\cal S}(K))$ 
\bit 
\itemI $T \subset f^{-1}(f(T)) = TK$ \ \ $(\upfa T \in {\cal S}(G))$ 
\itemII $\nu_H(T, f^{-1}(f(T))) \leq \nu_H(T, K) + 1$ \hsh $\nu_H(f^{-1}(f(T)), T) = 1$ 
\eit 
\eenum 
\eprop

\bnotb Suppose $f : G \to H$ is a group epimorphism and $K := {\rm Ker}\,f$.
Consider the following conditions for $T \in {\cal S}(G)$ : \\
\hsppp (i) \ $Q_m$ : $K \subset T^{\leq m}$ \hsp $(m \in \IZ_{\geq 0})$ \hsp (ii) \ $R$ : $T \cap K \in {\cal S}_g(K)$ \\
We have the subsets ${\cal S}_{Q_m}(G), \,{\cal S}_R(G) \subset {\cal S}(G)$. 
\enot

\bcorb Suppose $f : G \to H$ is a group epimorphism and $K := {\rm Ker}\,f$.
\benum[(1)] 
\item[{[I]}\ ] If $f$ is an isomorphism, then \\
\hsh $f_\ast : ({\cal S}''(G)^\ast, \widehat{\nu}_H) \lra ({\cal S}''(H)^\ast, \widehat{\nu}_H)$ is an isometric bijection
 and $f^\ast = f_\ast{}^{-1}$. 
\item[{[II]}\,]  
In each of the following cases, $f_\ast$ is a 1-Lipschitz surjection and 
$f^\ast$ is an isometric section of $f_\ast$. 
Moreover, $f_\ast$ is a quasi-isometric map and $f^\ast$ is a quasi-isometric inverse of $f_\ast$ 
\bit 
\item[(1)] $({\cal S}''_{Q_m}(G)^\ast - {\cal S}(K),  \widehat{\nu}_H) 
\bary[c]{c}
\mbox{\small $f_\ast$} \\[-1mm] 
\smash{\raisebox{0.5mm}{\makebox(20,0){\rightarrowfill}}} \\[-2mm] 
\smash{\raisebox{1mm}{\makebox(20,0){\leftarrowfill}}} \\[-1mm] 
\mbox{\small $f^\ast$}
\eary  ({\cal S}''(H)^\ast, \widehat{\nu}_H)$ \hsp 
$\circ$ \ \btab[t]{l}
$\widehat{\nu}_H(T, f^\ast f_\ast(T)) \leq \nu_H(T, K) + 1 \leq m + 1$ \\[2mm]
\hsp \hsh $(\upfa T \in {\cal S}''_{Q_m}(G)^\ast - {\cal S}(K))$
\etab 

\item[(2)] Suppose $K$ is a finite group. 
\bit 
\itemI $({\cal S}''_R(G)^\ast - {\cal S}(K),  \widehat{\nu}_H) 
\bary[c]{c}
\mbox{\small $f_\ast$} \\[-1mm] 
\smash{\raisebox{0.5mm}{\makebox(20,0){\rightarrowfill}}} \\[-2mm] 
\smash{\raisebox{1mm}{\makebox(20,0){\leftarrowfill}}} \\[-1mm] 
\mbox{\small $f^\ast$}
\eary  ({\cal S}''(H)^\ast, \widehat{\nu}_H)$
\vskip 1mm 
\bit 
\itema Since ${\cal S}(K)$ is a finite set,  there exists $m \in \IZ_{\geq 0}$ with $K = U^{\leq m}$ \ $(\upfa U \in {\cal S}_g(K))$. 
\itemb $\upfa T \in {\cal S}_R(G)$ \hsh $T \cap K \in {\cal S}_g(K)$ \hsh \tf \ $K = (T \cap K)^{\leq m} \subset T^{\leq m}$ \\[0.5mm] 
\tf \ ${\cal S}''_R(G)^\ast \subset {\cal S}''_{Q_m}(G)^\ast$.
\eit 
\vspace*{-2mm} 
\itemII Suppose $G$ is finitely generated. \hsh 
$({\cal S}''_{fg, R}(G)^\ast,  \widehat{\nu}_H) 
\bary[c]{c}
\mbox{\small $f_\ast$} \\[-1mm] 
\smash{\raisebox{0.5mm}{\makebox(20,0){\rightarrowfill}}} \\[-2mm] 
\smash{\raisebox{1mm}{\makebox(20,0){\leftarrowfill}}} \\[-1mm] 
\mbox{\small $f^\ast$}
\eary  ({\cal S}''_{fg}(H)^\ast, \widehat{\nu}_H)$ 
\vskip 1mm 
\bit 
\itema If $S \in {\cal S}''_{fg}(H)^\ast$, then 
$(f^{-1}(S))^\infty = f^{-1}(S^\infty) = f^{-1}(H) = G$. \hsp (cf. Fact~\ref{fact_f+f^{-1}}\,[II]).
\eit 
\eit 
\eit 
\eenum 
\ \ $\circ$ \ In the statements [I] and [II]\,(1), (2) we can replace $\widehat{\nu}_H$ by $\widehat{\rho}$. 
\ecor 

\subsection{The case of ${\cal S}_{nfg}(G)$} \mbox{} 

We extend the arguments in the previous subsection to the case of the map 
$f_\ast : ({\cal S}''_{nfg}(G)^\ast, \widehat{\rho}) \lra ({\cal S}''_{nfg}(H)^\ast, \widehat{\rho})$. 
Throughout this subsection we assume that  $f : G \to H$ is a group epimorphism between nontrivial groups and $K \equiv {\rm Ker}\,f$. 

\bnotb For $S \in {\cal S}(H)$ we say that $R \in {\cal S}(G)$ is a lift of $S$ with respect to $f$ if \\
\hsp \hsh $f(R) = S$ \ and \ $|R \cap f^{-1}(x)| = 1$ \ $(\upfa x \in S)$. \\
We use the following notations. 

\bit 
\itemI ${\cal S}(G; S, f) := \{ R \in {\cal S}(G) \mid \mbox{$R$ is a lift of $S$ with respect to $f$.} \}$
\itemII $R_x := R \cap f^{-1}(x) \in G$ \ \ $(R \in {\cal S}(G; S, f), \ x \in S)$
\eit 
\enot 

\bnotb 
\benum 
\item For $S\in {\cal S}(H)$ and $R \in {\cal S}(G; S, f)$ we have the map \\ 
\hsp \hsh $h = h_S^R : f^{-1}(S) \lra K$ : $h(x) = (R_{f(x)})^{-1}x$ 
\bit 
\itemI $x = R_{f(x)}h(x)$ \hsh $(\upfa x \in f^{-1}(S))$ 
\itemII $h(xa) = h(x)a$ \hsh $(\upfa x \in f^{-1}(S), a \in K)$
\itemiii $e \in S$, \ $R_e = e$ \LRA $h|_K = \id$, \ $h^{-1}(e) = R$
\eit 

\item There exists $\theta : {\cal S}(H) \lra {\cal S}(G)$ \ such that \ $\theta(S) \in  {\cal S}(G; S, f)$ \ for each $S \in {\cal S}(H)$. 
\bit 
\itemI If $g : H \to G$ is a (non-homomorphic) section of $f$, then we have the map \\
\hsp \hsh $\theta_g : {\cal S}(H) \lra {\cal S}(G)$ : $\theta_g(S) := g(S) \in {\cal S}(G; S, f)$ 
\eit 

\item Consider the following maps : \\
\hsh $\chi : {\cal S}(G) \lra {\cal S}(H \times K)$ : $\chi(T) := \{ (f(x), h_S(x)) \mid x \in T \} \subset f(T) \times K$, \\
\hsh $\omega : {\cal S}(H \times K) \lra {\cal S}(G)$ : $\omega(W) = \{ \theta(p(W))_uv \mid (u, v) \in W \}$, \\
\hsp where \ $S := f(T) \in {\cal S}(H)$ \ and \ $p : H \times K \lra H$ is the projection. 

\bit 
\itemI $T = \{ \theta(S)_u v \mid (u,v) \in \chi(T) \} \subset f^{-1}(S)$ \ \ $(T \in {\cal S}(G))$.  
\itemII $\omega \chi = \id$, $\chi \omega = \id$ \hsp ($\omega \chi (T) = T$, $\chi \omega(W) = W$) 
\eit 
\eenum 
\enot

\bfact\label{fact_theta} \mbox{}  
\benum
\item Let $S \in {\cal S}''(H)$, $R \in {\cal S}''(G; S, f)$ and $L, U \in {\cal S}''(K)$. Then, the following holds. 
\bit 
\itemI $f(R^L \cup U) = f(R^LU) = S$ \hsp (ii) \ $(R^L \cup U) \cap K =  (R^LU) \cap K = U$
\eit 
\vskip 1mm 
\item There exists a map \ $\theta : {\cal S}(H) \to {\cal S}(G)$ \ such that \\
\hsp \hsh $(\ast)$ \ $\theta(S) \in  {\cal S}(G; S, f)$ \ $(S \in {\cal S}(H))$ \ and \ $\theta({\cal S}''(H)) \subset {\cal S}''(G)$. \\
The following holds. 
\bit 
\itemI $\theta({\cal S}''(H)^\ast) \subset {\cal S}''(G)^\ast - {\cal S}(K)$ \hsp 
\itemII For each $S\in {\cal S}(H)$ we have \ $\theta(S) \in {\cal S}(G; S, f)$ \ and the map \\ 
\hspp $h_S := h_S^{\theta(S)} : f^{-1}(S) \to K$ : $h_S(x) = (\theta(S)_{f(x)})^{-1}x$.  
\itemiii $x = \theta(S)_{f(x)}h_S(x)$ \ $(\upfa x \in f^{-1}(S))$ 
\itemiv There exists a (non-homomorphic) section $g : H \to G$ of $f$ with $g(x^{-1}) = g(x)^{-1}$ $(x \in H)$ and $g(e) = e$.
Then, the map $\theta_g$ satisfies the above conditions $(\ast)$ and has the following property \\
 \hsp For $S \in {\cal S}(H)$ \hsh (a) $g(S^{-1}) = g(S)^{-1}$ \hsh (b) \ $S \in {\cal S}_s(H)$ \LRA  $g(S) \in {\cal S}_s(G)$ 
\eit 

\item Consider the following maps \\
$\phi : {\cal S}''(G) \lra {\cal S}''(H) \times {\cal S}''(K)$ : $\phi(T) = (f(T), T \cap K)$ \\
$\psi_L, \psi'_L : {\cal S}''(H)^\ast \times {\cal S}''(K)^\ast \lra {\cal S}''(G)^\ast  - {\cal S}(K)$ : \hsppp ($L \in {\cal S}''(K)$) \\
\hsppp $\psi_L(S, U) = \theta(S)^L \cup U, \ \ \psi'_L(S, U) = \theta(S)^LU$ 
\vskip 1mm 
\bit 
\itemI $\phi \psi_L = \phi \psi'_L = \id$ \hspp $\circ$ \ 
$f(R^L \cup U) = f(R^LU) = S$ \hsh $(R^L \cup U) \cap K =  (R^LU) \cap K = U$
\vskip 1mm 
\itemII $\widehat{\nu}_H(\psi_L(S,U), \psi'_L(S, U)) \leq 2$ \hspp 
$\circ$ \ $\theta(S)^L \cup U \subset \theta(S)^L U \subset (\theta(S)^L \cup U)^2$  
\eit 
\vskip 1mm 
\item We have the following restrictions of $\chi$ and $\omega$. \\
\hsp 
$\bary[t]{@{}c@{ \ \ }c@{ \ \ }c@{ \ \ }c@{ \ \ }cll@{}}
{\cal S}(G) & \supset & {\cal S}''(G) & \supset & {\cal S}''(G)  - {\cal S}(K) & & \\[-3mm] 
\rotatebox{-90}{$\lra$} \rotatebox{-90}{$\lla$} 
& \hspace*{-25mm} \smash{\raisebox{-4mm}{$\chi$}} \hspace*{5mm} \smash{\raisebox{-4mm}{$\omega$}} 
& \rotatebox{-90}{$\lra$} \rotatebox{-90}{$\lla$} 
& \hspace*{-26.5mm} \smash{\raisebox{-4mm}{$\chi$}} \hspace*{5mm} \smash{\raisebox{-4mm}{$\omega$}} 
& \rotatebox{-90}{$\lra$} \rotatebox{-90}{$\lla$} 
& \hspace*{-28mm} \smash{\raisebox{-4mm}{$\chi$}} \hspace*{5mm} \smash{\raisebox{-4mm}{$\omega$}} &
\smash{\raisebox{-4mm}{$\omega = \chi^{-1}$}} \\[6mm]
 {\cal S}(H \times K) & \supset & {\cal S}''(H \times K) & \supset & {\cal S}''(H \times K)  - {\cal S}(e \times K) & & 
\eary$ 
\vskip 2mm 
\eenum 
\efact 

\bnotb For a group $G$ we use the following notations and terminologies : 
\benum 
\item 
\bit 
\itemI ${\rm rank}_n\,G := \min \{ n \in \widetilde{\IZ}_{\geq 0} \mid \upexi A \subset G, \ |A| = n, \ G = \langle \!\langle A \rangle \!\rangle = (C_A)^\infty \}
\in \widetilde{\IZ}_{\geq 0}$ 
\itemII ${\rm diam}_{nfg} G := \min\,\{ \nu_H(S, G) \mid S \in {\cal S}_{nfg}(G) \} \in \widetilde{\IZ}_{\geq 0}$ 
\eit 
\item[] $\circ$ \ If $G$ is a nontrivial group, then these values are in $\widetilde{\IZ}_{\geq 1}$. 

\item $G$ is finitely normally generated \LLRA ${\cal S}_{nfg}(G) \neq \emptyset$ \LLRA ${\rm rank}_n\,G < \infty$ \\
\hspace*{53.5mm} \LLRA There exists $A \in {\cal S}_f(G)$ with  
$G = \langle \!\langle A \rangle \!\rangle = (C_A)^\infty$.  
\item $G$ is weak-uniformly finitely normally generated \\
\hsh \LLRA $\nu_H(S, G) < \infty$ for some/any $S \in {\cal S}_{nfg}(G)$ \LLRA ${\rm diam}_{nfg} G < \infty$ \\
\hsh \LLRA There exists $A \in {\cal S}_f(G)$ with $\nu_H(C_A,G) < \infty$ \ 
(i.e., $G = (C_A)^{\leq n}$ for some $n \in \IZ_{\geq 0}$). 
\bit 
\itemC If $S, T \in {\cal S}_{nfg}''(G)^\ast$, then $\nu_H(T, G) < \infty$ and $\nu_H(T, G) \leq \nu_H(T, S) \cdot \nu_H(S, G)$ \\
\hfill (cf. Example~\ref{exp_finite}\,(2)). 
\eit 

\item In \cite[Section 2]{KLM1} the following notations are introduced. 
\bit 
\itemI $\Delta(G) := \sup\,\{ \nu_H(S, G) \mid S \in {\cal S}_{nfg}(G) \} \in \widetilde{\IZ}_{\geq 0}$ 
\itemII $G$ is uniformly bounded. \LLRA $\Delta(G) < \infty$ 
%\LLRA There exists $n \in \IZ_{\geq 1}$ with $\nu_H(S, G) \leq n$ for any $S \in {\cal S}_{nfg}(G)$.  
\eit 
\eenum
\enot 

\bthm\label{prop_nfg} Suppose $f : G \to H$ is a group epimorphism between nontrivial groups and 
$K \equiv {\rm Ker}\,f$ is normally finitely generated. 
\bit 
\item[{[I]}\,] The map \ $f_\ast : ({\cal S}''_{nfg}(G)^\ast, \widehat{\rho}) \lra ({\cal S}''_{nfg}(H)^\ast, \widehat{\rho})$ \  is a 1-Lipschitz surjection with a section.  
\item[{[II]}] If $K$ is weak-uniformly normally finitely generated and $m_0 := {\rm diam}_{nfg}\,K \in \IZ_{\geq 0}$, then \\
the following holds. 
\benum 
\item[1)] The map $f_\ast : {\cal S}''_{nfg}(G)^\ast \lra  {\cal S}''_{nfg}(H)^\ast$ 
has a section $\eta : {\cal S}''_{nfg}(H)^\ast \lra {\cal S}''_{nfg, Q_{m_0}}(G)^\ast$.   
\item[2)] For any $m \in \IZ_{\geq m_0}$ the map \ \ 
$f_\ast : ({\cal S}''_{nfg, Q_m}(G)^\ast, \widehat{\rho}) \lra ({\cal S}''_{nfg}(H)^\ast, \widehat{\rho})$ \ \ 
is a quasi-isometry and $\eta$ is a quasi-isometric inverse of $f_\ast$. 
\bit 
\itemI $\widehat{\rho}(\eta \times \eta) \leq \widehat{\rho} + \log (m_0+1)$ \hsp 
(ii) \ $\widehat{\rho}(\eta f_\ast, \id) \leq \log(m+1)$
\eit 
\eenum 
\eit 
\ethm 

\bpfb 
The assertions follow from the following observations. 
\benum[(1)]
\item[{[I]}\,]  
\bit 
\item[(1)] (i) \ $f_\ast({\cal S}''_{nfg}(G)^\ast) \subset {\cal S}''_{nfg}(H)^\ast$ \ (cf.~Fact~\ref{fact_f+f^{-1}}\,[I](3)) \hsf 
(ii) \ $f_\ast$ is 1-Lipschitz \ (cf.~Proposition~\ref{prop_f^ast+f_ast}\,(1)). 
\item[(2)] Since $K$ is normally finitely generated, there exists $B \in {\cal S}_{f, s}''(K)$ with $K = (C_B)^\infty = (B^K)^\infty$ in $K$. 
\bit 
\itemI Let $U := C_B = B^G$ in $G$. Then $U \in {\cal S}''_{s, fc}(G) \cap {\cal S}(K)$ and $K = U^\infty$.
\eit 
\item[(3)] We define a map  
\hsh $\eta : {\cal S}''_{nfg}(H)^\ast \lra {\cal S}''_{nfg}(G)^\ast$ \hsh by the following rule :
\vskip 1mm 
\bit 
\item[$(\sharp)$] \ For each $S \in {\cal S}''_{nfg}(H)^\ast$ \ 
\btab[t]{cl}
(a) & choose $A \in {\cal S}''_{f,s}(H)$ with $S = C_A = A^H$ \ and \\[2mm] 
(b) & with the lift $\theta(A) \in {\cal S}''_{f}(G, A, f)$ \ we set \ \ $\eta(S) := C_{\theta(A)} \cup U$. 
\etab  
\eit 
\vskip 1mm 
\bit 
\itemI $f(C_{\theta(A)}) = C_{f(\theta(A))}) = C_A = S$. \hsh (cf. Fact~\ref{fact_f+f^{-1}}\,[I]) 
\eit 
\vskip 1mm 
\item[(4)] $T := \eta(S) \in {\cal S}''_{nfg}(G)^\ast$
\vskip 1mm 
\bit 
\itemI $T \in {\cal S}_g(G)$ : 
\hsf  
\btab[t]{l}
$f((C_{\theta(A)})^\infty) = f(C_{\theta(A)})^\infty = S^\infty = H$. \hsh 
(cf. Fact~\ref{fact_f+f^{-1}}\,[I]) \\[2mm]
\tf \ $G = (C_{\theta(A)})^\infty K \subset (C_{\theta(A)})^\infty U^\infty \subset (C_{\theta(A)} \cup U)^\infty = T^\infty$.
\etab 
\vskip 1.5mm 
\itemII $T \in {\cal S}_{s, fc}(G)$ : \hsh $T = C_{\theta(A)} \cup U = C_{\theta(A)} \cup C_B = C_{\theta(A) \cup B}$ \ in $G$ \hsh 
$\theta(A) \cup B \in {\cal S}_{f}(G)$  
\eit 
\item[(5)] $f_\ast \eta = \id$ : \hsh 
$f(\eta(S)) = f(C_{\theta(A)} \cup U) = f(C_{\theta(A)}) \cup f(U) = S \cup \{ e \} = S$. 
\eit 
\vskip 2mm 

\item[{[II]}] 
\bit 
\item[(1)] By the assumption there exists $\upexi B \in {\cal S}_{f, s}''(K)$ with $K = (C_B)^{\leq m_0} = (B^K)^{\leq m_0}$ in $K$. 
\bit 
\itemI Let $U := C_B = B^G$ in $G$. Then, \ $U \in {\cal S}''_{s, fc}(G)^\ast \cap {\cal S}(K)$ \ and \ $K = U^{\leq m_0}$. % = U^{\leq m}$. 
\eit 
\vskip 1mm 
\item[(2)] $T, T' \in {\cal S}''_{Q_{n}}(G)^\ast - {\cal S}(K)$ $(n \in \IZ_{\geq 0})$,\ $S = f(T), S' = f(T')$ \LRA 
$\nu_H(T, T') \leq (n+1) \nu_H(S, S')$. 
\bit 
\itemI $K \subset T^{\leq n}$ \hsh $T \subset f^{-1}(S) = TK \subset T^{\leq n+1}$ \\
\tf \ $\nu_H(T, f^{-1}(S)) \leq n+1$, \ $\nu_H(f^{-1}(S), T) = 1$ \hsh 
\tf \ $\widehat{\nu}_H(T, f^{-1}(S)) \leq n+1$
\item[] $\circ$ \ Similar for $T'$ and $S'$. 
\itemII $\nu_H(f^{-1}(S), f^{-1}(S')) = \nu_H(S, S')$ \hsh (cf.~Proposition~\ref{prop_f^ast+f_ast}\,(2)) 
\itemiii $\nu_H(T, T') \leq \nu_H(T, f^{-1}(S)) \cdot \nu_H(f^{-1}(S), f^{-1}(S'))  \cdot \nu_H(f^{-1}(S'), T') \leq (n+1) \nu_H(S, S')$.  
\eit 
\item[(3)] 
We define the map  
\hsf $\eta : {\cal S}''_{nfg}(H)^\ast \lra {\cal S}''_{nfg}(G)^\ast$ \hsf as in [I]. 
\vskip 1mm 
\bit 
\itemI $\eta(S) \in {\cal S}''_{nfg, Q_{m_0}}(G)^\ast$, \ since \ $K = U^{\leq m_0} \subset \eta(S)^{\leq m_0}$. 
\hsp (ii) \ $f(\eta(S)) = S$. 
\itemiii $\nu_H(\eta(S), \eta(S')) \leq (m_0+1) \nu_H(S, S')$ \ \ $(S, S' \in {\cal S}''_{nfg}(H)^\ast)$ \hsh (by (2))
\itemiv $\widehat{\nu}_H(\eta \times \eta) \leq (m_0+1) \widehat{\nu}_H$ \hsp \tf \ 
$\widehat{\rho}(\eta \times \eta) \leq \widehat{\rho} + \log (m_0+1)$ \eit 
\vskip 1mm 
\item[(4)] $f_\ast : ({\cal S}''_{nfg, Q_m}(G)^\ast, \widehat{\nu}_H) \lra 
({\cal S}''_{nfg}(H)^\ast, \widehat{\nu}_H)$  is a 1-Lipschitz surjection with a section $\eta$. 
\vskip 1mm 
\bit 
\itemI $\eta : ({\cal S}''_{nfg}(H)^\ast, \widehat{\nu}_H) \lra ({\cal S}''_{nfg, Q_m}(G)^\ast, \widehat{\nu}_H)$ is a quasi-isometric map. 
\hsh (ii) \ $f_\ast \eta = \id$
\itemiii For $T \in {\cal S}''_{nfg, Q_m}(G)^\ast$ and $S := f(T)$ \\
\hsp $\widehat{\nu}_H(T, \eta(S)) \leq m + 1$ \ \ by (2) \hsh \tf \ $\widehat{\rho}_H(T, \eta(S)) \leq \log(m + 1)$.
\eit 
\eit 
\eenum 
\vspace*{-6mm} 
\epf 

\subsection{Semidirect products} \mbox{} 

This subsection is devoted to the case of semidirect products (or the case of group epimorphisms with homomorphic sections). 
Suppose a group $G$ is the semidirect product of a subgroup $H$ and a normal subgroup $K$ ($G = H \ltimes K$). 
Let $p : G \to H$ and $q : G \to K$ denote the projections onto the factors $H$ and $K$.
The map $p$ is a group epimorphism with $p|_H = \id_H$ and ${\rm Ker}\,p = K$.
Any $S \in {\cal S}(H)$ is its canonical lift with respect to $p$, to which we can apply the arguments in the previous subsection.
We set \\
\hspp ${\cal S}_P''(K)^G := {\cal S}_P''(K) \cap {\cal S}_{c}(G) \subset {\cal S}_P''(K)$ \hsh (cf. Notation~\ref{not_S_p}). 

\bfactb\label{fact_S^KcupU} 
\benum 
\item Let $S \in {\cal S}''(H)$, $U \in {\cal S}''(K)$ and $L \in {\cal S}''(K)$.
\bit 
\itemI $S \subset S^L \subset S^L \cup U \subset S^L U$ 
\itemII $p(S^L) = p(S^L \cup U) = p(S^L U) = S$ \hsh $(S^L \cup U)\cap K = (S^LU)\cap K = U$ \\
$S^L \cap H = S$ \hsh $S^LU \cap H = S$ 
\itemiii $S^KU = US^K$ \bec \ $S^K$ : $K$-conjugate-invariant \hsh $U \subset K$ 
\item[] $\circ$ \ $p(S^L) = p(S)^{p(L)} = S^e = S$, \hsh $p(S^LU) = p(S^L)p(U) = S e = S$ 
\item[] $\circ$ \ $S^K \cap K = \{ e \}$, \ \ 
\item[] $\circ$ \ $S^K = ((S - \{ e \}) \cup \{ e \})^K = (S - \{ e \})^K \cup \{ e \}^K 
= (S - \{ e \})^K \cup \{ e \} \subset (S - \{ e \})K \cup \{ e \}$ 
\item[] $\circ$ \ $S^KU \subset ((S - \{ e \})K \cup \{ e \})U \subset (S - \{ e \})K \cup U$ 
\itemiii $AK \cap BK = (A \cap B)K$ \ \ $(A, B \subset H)$, \hsh $(H^\times K) \cap K = \emptyset$, \hsh $A^K \subset A^KK\subset AK$ 
\ \ $(A \subset G)$ 
\eit 

\item Let $S \in {\cal S}''_{s,c}(H)$ and $U \in {\cal S}''_{s,c}(K)^G := {\cal S}_{s, c}''(K) \cap {\cal S}_{c}(G)$
\bit 
\itemI $S^K \in {\cal S}''_{s,c}(G)$ \hsp (ii) \ $S^K \cup U, \ S^KU \in {\cal S}''_{s,c}(G)$ \hsp (cf. Fact~\ref{fact_formula}\,(6)(vi))
\item[] $\circ$ \ $(S^K)^G = S^{KG} = S^G = S^{HK} = (S^H)^K = S^K$ \hsh $(S^K)^{-1} = (S^{-1})^K = S^K$
\eit 
\item Let $S \in {\cal S}''_{g}(H)$ and $U \in {\cal S}''_{g}(K)$. 
\bit 
\itemI $S \cup U, \,SU \in {\cal S}''_{g}(G)$
\item[] $\circ$ \ 
$(S \cup U)^\infty \supset S^\infty \cup U^\infty = H \cup K$, \ \ 
$(S \cup U)^\infty = ((S \cup U)^\infty)^2 \supset (H \cup K)^2 \supset HK = G$ 

\itemII $S \in {\cal S}''_{fg}(H)$, \ $U \in {\cal S}''_{fg}(K)$ \LRA $S \cup U \in {\cal S}''_{fg}(G)$, \ $SU \in {\cal S}''_{f, \,g}(G)$
\vskip 1mm 
\itemiii $S \in {\cal S}''_{nfg}(H)$, \ $U \in {\cal S}''_{g}(K)$, \ $U = B^G$ \ \ $(\upexi B \in {\cal S}_{f,s}(K))$ \LRA $S^K \cup U \in {\cal S}''_{nfg}(G)$
\item[] $\circ$ \ $S = A^H$ \ \ $(\upexi A \in {\cal S}_{f,s}(H))$ \hsh 
$U = B^G$ \ \ $(\upexi B \in {\cal S}_{f,s}(K))$ \hsh $e \in S$, $e \in U$ \hsh \tf \ $e \in A$, \ $e \in B$ 
\item[] $\circ$ \ $S^K = (A^H)^K = A^{HK} = A^G$ \ \ \ $S^K \cup U = A^G \cup B^G = (A \cup B)^G = C_{A \cup B}$ \ \ \ 
$A \cup B \in {\cal S}_{f,s}(G)$
\item[] $\circ$ \ $S^KU \supset S^K \cup U \supset S \cup U \in {\cal S}''_{g}(G)$ 
\item[] $\circ$ \ $S^KU = A^GB^G = B^GA^G\supset (AB \cup BA)^G \supset A^G \cup B^G = S^K \cup U$ 
\eit 
\eenum 
\efact 

\bfactb Let $T \in {\cal S}''(G)$. We have $p(T) \in {\cal S}''(H)$ and $T \cap K \in {\cal S}''(K)$. 
\benum 
\item The following implication holds for the conditions  $P = f, \ g,\ s, \ c, \ fc, \ fg, \ ng, \ nfg$ : \\
\hsh $T \in {\cal S}''_P(G)$ \LRA $p(T) \in {\cal S}''_P(H)$ \hspp (cf. Fact~\ref{fact_f+f^{-1}}\,[I])
\item The following implication holds for the conditions  $P = f, \ s, \ c$ : \\
\hsh $T \in {\cal S}_P''(G)$ \LRA $T \cap K \in {\cal S}_P''(G)$ \hspp (cf. Fact~\ref{fact_formula}\,(5), (6))
\eenum 
\efact 

\bnotb Consider the following maps. 
\benum 
\item $\phi : {\cal S}''(G) \lra {\cal S}''(H) \times {\cal S}''(K)$ : $\phi(T) = (p(T), T \cap K)$
\item 
$\bary[t]{@{}c@{ \ }c@{ \ }c@{ \ }c@{ \ }c@{}}
\psi = \psi_L : & {\cal S}''(H)^\ast \times {\cal S}''(K)^\ast & \lra & {\cal S}''(G)^\ast : & \psi(S, U) = S^L \cup U \\[1.5mm] 
\psi' = \psi'_L : & {\cal S}''(H)^\ast \times {\cal S}''(K)^\ast & \lra & {\cal S}''(G)^\ast : & \psi'(S, U) = S^LU
\eary$ \hsp ($L \in {\cal S}''(K)$)
\vskip 2mm 

\bit 
\itemI $\phi \psi = \phi \psi' = \id$ \hsh (cf. Fact~\ref{fact_S^KcupU}\,(1)) 
\itemII $\widehat{\nu}_H(\psi(S,U), \psi'(S, U)) \leq 2$ \\
\ \ $\circ$ \ $S^L \cup U \subset S^L U \subset (S^L \cup U)^2$ \hsh 
\tf \ $\nu_H(S^L \cup U, S^L U) \leq 2$, \ $\nu_H(S^L U, S^L \cup U) \leq 1$
\itemiii When $G$ is the product of $H$ and $K$, we have $S^L = S$. 
\eit 
\item We use the following notations : 
\bit 
\itemI $\nu_H^G$, $\nu_H^H$, $\nu_H^K$ : the multiplicative as metric $\nu_H$ on 
${\cal S}''(G)^\ast$, ${\cal S}''(H)^\ast$, ${\cal S}''(K)^\ast$ respectively. 
\itemII $\nu_H^{H, K} \in {\cal M}_{m,as}({\cal S}''(H)^\ast \times {\cal S}''(K)^\ast)$ : the product metric of $\nu_H^H$ and $\nu_H^K$ \\
\hsp $\nu_H^{H, K}((S,U), (S', U')) = \max\{ \nu_H^H(S,S'), \nu_H^K(U,U') \}$ 
\itemiii $\mu, \mu' \in {\cal M}_{m,as}({\cal S}''(H)^\ast \times {\cal S}''(K)^\ast)$ : the metrics induced from $\nu_H^G$ 
by the injections $\psi$ and $\psi'$ 
\item[] \hsh $\bary[t]{@{}l@{}}
\mu((S,U), (S', U')) = \nu_H^G(\psi(S,U), \psi(S', U')) \\[1.5mm]
\mu'((S,U), (S', U')) = \nu_H^G(\psi'(S,U), \psi'(S', U'))
\eary$ \ \ $((S,U), (S', U') \in {\cal S}''(H)^\ast \times {\cal S}''(K)^\ast)$ 
\eit 
\eenum 
\enot

We are interested in the description of ${\rm Im}\,\psi$ and a deviation of $\mu$ and $\mu'$ from $\nu_H^{H, K}$. 

\bprop\label{prop_sdprod} Let $(S,U), (S', U') \in {\cal S}''(H)^\ast \times {\cal S}''(K)^\ast$ \ and \ $L \in {\cal S}''(K)$. \\
\ {[I]} the case of $\psi$ : \hsh Let $T = S^L \cup U$, $T' =(S')^L \cup U'$. 
\benum 
\item For $X \in {\cal S}''(G)$ \hsh $X \in {\rm Im}\,\psi$ \LLRA $X = (X \cap H)^L \cup (X \cap K)$, \ $X \cap H \in {\cal S}''(H)^\ast$, \ 
$X \cap K \in {\cal S}''(K)^\ast$
\bit 
\itemC $T \cap H = p(T) = S$, \ $T \cap K = U$ 
\eit 
\item 
\bit 
\itemI $\psi : {\cal S}''_P(H)^\ast \times {\cal S}''_P(K)^\ast \lra {\cal S}''_P(G)^\ast$ \ in the following cases : 
\item[] (a) \ $P = g$ \hsh (b) \ $P = s$ \hsh (c) \ $P = f$ and $L \in {\cal S}''_f(K)$ \hsh (d) \ $P = fg$ and $L \in {\cal S}''_f(K)$
\item[] \hspace*{-10mm} When $L = K$ : 
\itemII $\psi : {\cal S}''_P(H)^\ast \times ({\cal S}''_P(K)^\ast \cap {\cal S}_c(G)) \lra {\cal S}''_P(G)^\ast$ \ \ for \ (a) \ $P = c$, \ \ (b) \ $P = ng$ 
\item[] \ \ $\circ$ \ If \ $S \in {\cal S}''_c(H)$, then \ $(S^K)^G = S^{KG} = S^G = S^{HK} = (S^H)^K = S^K$. 
\itemiii $\psi : {\cal S}''_{fc}(H) \times ({\cal S}''(K) \cap {\cal S}''_{fc}(G)) \lra {\cal S}''_{fc}(G)$
\item[] \ \ $\circ$ \ If \ $S = A^H$ \ $(A \in {\cal S}_{f}(H))$, \ \ $U = B^G$ \ $(B \in {\cal S}_{f}(K))$, then \\
\hsh \ $S^K = (A^H)^K = A^{HK} = A^G$ \hsh $S^K \cup U = A^G \cup B^G = (A \cup B)^G$
\itemiii $\psi : {\cal S}''_{nfg}(H) \times ({\cal S}''_{s, g}(K) \cap {\cal S}_{fc}(G)) \lra {\cal S}''_{nfg}(G)$ 
\vskip 1mm 
\bit 
\itemC ${\rm Im}\,\psi \subset \{ T \in {\cal S}''_{nfg}(G) \mid T \cap K \in {\cal S}''_{s, g}(K) \cap {\cal S}_{fc}(G)\}$  
\eit 
\eit 
\vskip 1mm 
\item $\mu((S,U), (S', U')) = \nu_H^G(T, T') = \max \{ \nu_H^H(S, S'), \nu_H^G(T, U') \}$ 
\bit 
\itemI For $n \in \IZ_{\geq 0}$ \hsh 
$\nu_H^G(T, T') \leq n$ \LLRA $T' \subset T^{\leq n}$ \LLRA [$S' \subset S^{\leq n}$, \ $U' \subset T^{\leq n}$] \\
\hspace*{50mm} \LLRA [$\nu_H^H(S, S') \leq n$, \ $\nu_H^G(T, U') \leq n$] 
\vskip 1mm 
\item[] \ \ $\circ$ \ $T' \subset T^{\leq n}$ \LLRA [$S' \subset S^{\leq n}$, \ $U' \subset T^{\leq n}$]
\bit 
\item[] \hspace*{-2mm} $(\Lra)$ \ $S' = p(T') \subset p(T^{\leq n}) = p(T)^{\leq n} = S^{\leq n}$ 
\item[] \hspace*{-2mm} $(\Lla)$ \ $(S')^L \subset (S^{\leq n})^L \subset (S^L)^{\leq n} \subset T^{\leq n}$ 
\eit 
\eit 
\eenum
\vskip 1mm 
\ {[II]} the case of $\psi'$ : \hsh Let $T = S^LU$, $T' =(S')^LU'$. 
\benum 
\item For $X \in {\cal S}''(G)$ \hsh $X \in {\rm Im}\,\psi'$ \LLRA $X = (X \cap H)^L(X \cap K)$, \ $X \cap H \in {\cal S}''(H)^\ast$, \ 
$X \cap K \in {\cal S}''(K)^\ast$
\bit 
\itemC $T \cap H = p(T) = S$, \ $T \cap K = U$ 
\eit 
\item 
\bit 
\itemI $\psi' : {\cal S}''_P(H)^\ast \times {\cal S}''_P(K)^\ast \lra {\cal S}''_P(G)^\ast$ \ in the following cases : 
\item[] (a) \ $P = g$ \hsh (b) \ $P = s$ \hsh (c) \ $P = f$ and $L \in {\cal S}''_f(K)$ \hsh (d) \ $P = fg$ and $L \in {\cal S}''_f(K)$
\bit 
\item[\bec] (b) \ If $S^{-1} = S$ and $U^{-1} = U$, then \\
\hsp $(S^KU)^{-1} = U^{-1}(S^K)^{-1} = U(S^{-1})^K = US^K = S^KU$ 
\hsh (cf. Fact~\ref{fact_formula}\,(1)(vi)) 
\eit 

\item[] \hspace*{-10mm} When $L = K$ : 
\itemII $\psi' : {\cal S}''_P(H)^\ast \times ({\cal S}''_P(K)^\ast \cap {\cal S}_c(G)) \lra {\cal S}''_P(G)^\ast$ \ \ 
for \ (a) \ $P = c$, \ \ (b) \ $P = ng$
\bit 
\item[\bec] (a) \ If $S^H = S$ and $U^G = U$, then $(S^K)^G = S^K$ and $(S^KU)^G \subset (S^K)^GU^G = S^KU$. 
\eit 
\eit 

\item $\mu'((S,U), (S', U')) = \nu_H^G(T, T')$ : \\[0.5mm] 
\hsp \hsh $\max \{ \nu_H^H(S, S'), \nu_H^G(T, U') \} \leq \nu_H^G(T, T') \leq 2\max \{ \nu_H^H(S, S'), \nu_H^G(T, U') \}$
\vskip 0.5mm 
\bit 
\itemI For $n \in \IZ_{\geq 0}$ \hsh 
$T' \subset T^{\leq n}$ $\relsb{}{\Lra}{(\sharp)}$ [$S' \subset S^{\leq n}$, \ $U' \subset T^{\leq n}$] $\relsb{}{\Lra}{(\flat)}$ $T' \subset T^{\leq 2n}$
\item[] \hsh $(\sharp)$ \ $S' = p(T') \subset p(T^{\leq n}) = p(T)^{\leq n} = S^{\leq n}$ 
\hsh $(\flat)$ \ $(S')^L \subset (S^{\leq n})^L \subset (S^L)^{\leq n} \subset T^{\leq n}$
\eit 
\eenum
\eprop 

As expected, we can see that $\mu = \mu' = \nu_H^{H, K}$, when $G$ is the direct product of the subgroups $H$ and $K$. 

\bfact\label{fact_HtimesK} In the case $G = H \times K$, the following holds. 
\benum 
\item $[H, K] = \{ e \}$ 
\item The map $q : G \to K$ is also a group epimorphism with $q|_K = \id$ and ${\rm Ker}\,q = H$. 
\item Let $S \in {\cal S}''(H)$ and $U \in {\cal S}''(K)$.
\bit 
\itemI $q(S \cup U) = q(SU) = U$ 
\itemII  $S^L = S$ \ $(\upfa L \in {\cal S}''(K))$ \ \ \ \tf \ $\psi$, $\psi'$ are independent of $L$. 
\itemiii $SU = US$ \hsh $(SU)^n = S^nU^n$, \hsh $(S \cup U)^n \supset S^n \cup U^n$ \hsh $(n \in \IZ_{\geq 0})$
\itemiv $S^K = S$ \hsh $S^G = S^{KH} = (S^K)^H = S^H$ 
\itemv $U^H = U$ \hsh $U^G = U^{HK} = (U^H)^K = U^K$ 
\item[(vi)] $(SU)^G = S^HU^K$ \hsh 
$\circ$ \ $(xy)^{ab} = x^{ab}y^{ab} = (x^b)^a(y^a)^b = x^ay^b$ \ \ $((x,a, y, b) \in S \times H \times U \times K)$
\eit 
\item 
\bit 
\itemI ${\cal S}_c''(H) = {\cal S}''(H) \cap {\cal S}_{c}(G)$ 
\itemII ${\cal S}_c''(K) = {\cal S}''(K) \cap {\cal S}_{c}(G)$ \hsh 
${\cal S}''_P(K) \cap {\cal S}_c(G) = {\cal S}''_{P,c}(K)$ 
\itemiii ${\cal S}''(K) \cap {\cal S}''_{fc}(G) = {\cal S}''_{fc}(K)$
\itemiv ${\cal S}''_{s, g}(K) \cap {\cal S}_{fc}(G) = {\cal S}''_{nfg}(K)$ 
\eit 
\eenum 
\efact 

\bprop\label{prop_prod} Suppose $G = H \times K$. 
For $(S,U), (S', U') \in {\cal S}''(H)^\ast \times {\cal S}''(K)^\ast$ \ we set \\
\hsp $T = S \cup U$, $T' =S' \cup U'$ \ in the case of $\psi$ \ \ and \ \ $T = SU$, $T' =S'U'$ \ in the case of $\psi'$. 
\benum 
\vskip 1mm 
\item  $\mu = \mu' = \nu_H^{H, K}$ \hsh (i.e., $\mu((S,U), (S', U')) = \nu_H^G(T, T') = \max \{ \nu_H^H(S, S'), \nu_H^K(U, U') \}$) 
\vskip 1mm 
\bit 
\itemI For $n \in \IZ_{\geq 0}$ \hsh $T' \subset T^n$ \LLRA $S' \subset S^n$, \ $U' \subset U^n$ 
\item[] \bec $(\Lra)$ \ $S' = q(T') \subset q(T^n) = q(T)^n = S^n$, \ \ $U' = q(T') \subset q(T^n) = q(T)^n = U^n$.
\item[] \hsh \ $(\Lla)$ \ $(S \cup U)^n \supset S^n \cup U^n \supset S' \cup U'$,  
\hsh $(SU)^n = S^nU^n \supset S'U'$ 
\eit 
\item[] $\circ$ \ The following maps are isometric injections : \\[1mm] 
\hsp $\psi, \psi' : ({\cal S}''(H)^\ast \times {\cal S}''(K)^\ast, \widehat{\rho}^{H, K}) \lra ({\cal S}''(G)^\ast, \widehat{\rho}^G)$ : 
\ \ $\psi(S, U) = S \cup U$, \ $\psi'(S, U) = SU$ \\[1mm] 
 \hsp Here, \ $\widehat{\rho}^{H, K} = \log \widehat{\nu}_H^{H, K} = \max \{ \widehat{\rho}^H, \widehat{\rho}^K \}$ 
\vskip 1mm 
\item $\psi, \psi' : {\cal S}''_P(H)^\ast \times {\cal S}''_P(K)^\ast \lra {\cal S}''_P(G)^\ast$ \hsh for $P = g, \ s, \ f, \ c, \ fc, \ fg, \ ng, \ nfg$. 
\bit 
\itemI The assertions follow from Fact~\ref{fact_HtimesK}\,(4) and Proposition~\ref{prop_sdprod}\,[I](2), [II](2) except the cases for 
$\psi'$ and $P = fc, \ nfg\, (= g + s + fc)$.  

\itemII the case of $\psi'$ and $P = fc$ : \ 
If \ $S = A^H$ \ $(A \in {\cal S}_{f}(H))$, \ $U = B^K$ \ $(B \in {\cal S}_{f}(K))$, then \ $SU = A^HB^K = (AB)^G$ \ 
(cf. Fact~\ref{fact_HtimesK}\,(3)(iv)).
\eit 
\eenum 
\eprop 

\section{Word metrics in group actions}

\subsection{Word metrics in group actions} \mbox{} 

Suppose $X$ is a set, $G$ is a group, $e$ is the identity element of $G$ and $\eta : X \times G \to X$ is a right group action. 

\bnotb We use the following notations. Let $x \in X$, $A \subset X$ and $a \in G$, $S \subset G$.  
\benum 
\item We have the maps \hsh $\eta_x : G \to xG \subset X$ : $\eta_x(a) = xa$ \ \ and \ \ $\eta^a : X \approx X$ : $\eta^a(x) = xa$. 

\bit
\itemC The map $\eta_x$ is injective if and only if the action of $G$ on the orbit $xG$ is free. 
\eit 
\item  $AS := \{ xs \mid x \in A, s \in S \} \subset X$
\eenum 
\enot

\bfactb Let $S, T, U, S_\lambda \subset G$ \ $(\lambda \in \Lambda)$ \ and \ $x,y \in X$, $A, B \subset X$. 
\benum 
\item (i) \ $y \in xS$ \LLRA $x \in yS^{-1}$ 
\hsh (ii) \ $y \in xS^n$ \LLRA $x \in y(S^{-1})^n = y(S^n)^{-1}$ \ \ $(n \in \widetilde{\IZ}_{\geq 0})$

\item 
(i) \ $A \subset B$, \ $S \subset T$ \LRA $AS \subset BT$ 
\hsh (ii) \ $(AS)T = A(ST)$ 
\hsh (iii) \ $A \big( \bigcup_{\lambda \in \Lambda} S_\lambda \big) = \bigcup_{\lambda \in \Lambda} AS_\lambda$
\item (i) \ $AS^0 = A$ \hsh (ii) \ $AS^\infty = \bigcup_{n = 0}^\infty A S^n$
\eenum 
\efact 

\bdefnb Let $S \subset G$.  
\benum 
\item the word metric on $X$ with respect to $S$ : \\[2mm] 
\hsp $d_S : X \times X \lra \widetilde{\IZ}_{\geq 0}$ : \ \  
$d_S(x,y) = \left\{ \hspace*{-1mm} \bary[c]{ll}
\min \{ n \in \IZ_{\geq 0} \mid y \in xS^n \} & (y \in xS^\infty) \\[1mm]
\infty & (y \in X - xS^\infty)
\eary \right.$ 
\vskip 2mm 

\item the Hausdorff metric on ${\cal S}(X)$ induced from $d_S$ : \\[1mm] 
\hsp $(d_S)_H : {\cal S}(X) \times {\cal S}(X) \lra \widetilde{\IZ}_{\geq 0}$ : \\[2mm]  
\hsp \hsh    
$(d_S)_H(A, A') := \left\{ \hspace*{-1mm} \bary[c]{ll}
\min \{ n \in \IZ_{\geq 0} \mid A' \subset AS^{\leq n} \} & (\upexi  n \in \IZ_{\geq 0} \ \ \mbox{s.t.} \ A' \subset AS^{\leq n}) \\[2mm]
\ \infty & (\mbox{otherwise}) 
\eary \right.$
\vskip 2mm 
\bit 
\itemI For $n \in \IZ_{\geq 0}$ \hsh $(d_S)_H(A, A') \leq n$ \LLRA $A' \subset A S^{\leq n}$ 
\eit 
\eenum
\edefn 

\bfact\label{fact_d_S} Let $x,y \in X$, $a \in G$, $S, T \subset G$ and $n \in \IZ_{\geq 0}$. 
\benum
\item $d_S \in {\cal M}_{as}(X)$, \ $d_S$ is nondegenerate \hsh \tf \ $\widehat{d_S} \in {\cal M}(X)$
\bit 
\itemI $d_S(x,y) \leq n$ \LLRA $y \in xS^{\leq n}$ \hsp (ii) \ $d_S(x,y) = 0$ \LLRA $x = y$  
\vskip 1mm 
\itemiii $B_{d_S}(x,n) := \{ y \in X \mid d_S(x,y) \leq n \} = xS^{\leq n}$ 
\eit 
\item \bit 
\itemI $d_S(x,y) = d_{S^{-1}}(y,x)$ \hsp (ii) \ $S = S^{-1}$ \LRA $d_S \in {\cal M}(X)$ 
\eit 
\item 
\bit 
\itemI $d_{S^a}(xa, ya) = d_S(x,y)$ \hsp (ii) \ $d_S(xa, ya) = d_S(x,y)$ \  if $S$ is conjugation-invariant.
\eit 
\item 
\bit 
\itemI $d_S = d_{S \cup \{ e \}}$ \hsp (ii) \ $d_S(x, \ )^{-1}(1) = xS - x$ 
\eit 
\item $S \subset T$ \LRA $d_S^X \geq d_T^X$. 
\eenum
\efact 

\bexpb 
\benum 
\item In the case $(X, \eta) = (G, R)$, where $R$ denotes the action of $G$ on itself by the right translation,  
for any $S \subset G$ the metric $d_S$ is just the word metric defined in Definition~\ref{defn_d_S}. 

\item Suppose $X = (X, \ast)$ is a set with a binary operation. 
Consider the group $G = {\rm Aut}(X)$ with the multiplication $f \cdot g = g \circ f$ \ $(f,g \in G)$. 
\bit 
\itemI The group $G$ has the canonical right action on $X$ : \ $x \cdot f = f(x)$ \ $(x \in X, f \in G)$. 
\itemII Any subset $S \subset G$ induces the word metric $d_S$ on $X$. 
\eit 
When $X = (X, \ast)$ is a quandle, it has the symmetries : \\
\hspp \hsh $\sigma_a \in {\rm Aut}(X)$ \ $(a \in X)$ : \ \ $\sigma_a(x) = x \ast a$ \ \ $(x \in X)$ \\
The group $G$ includes the following subgroups : \\
\hspp $G > 
{\rm Inn}(X) := \langle \sigma_a \ (a \in X) \rangle > 
{\rm Dis}(X) := \langle \sigma_a \sigma_b^{-1} \ (a,b \in X) \rangle$. \\
For example, we can take $S$ as the generating sets of these subgroups.  
\eenum
\eexp 

\bfact\label{fact_eta} 
Each $S \in {\cal S}(G)$ induces the word metrics $d_S \equiv d_S^G \in {\cal M}_{as}(G)$ and $d_S \equiv d_S^X \in {\cal M}_{as}(X)$.
\benum 
\item 
For $x \in X$ we have the map $\eta_x : (G, d_S) \,\lra\hspace{-5.5mm}\lra \,  (xG, d_S|)$. 
\vskip 0.5mm 
\bit 
\itemI $\eta_x$ is a 1-Lipschitz map. \hsp (ii) \ If $\eta_x$ is injective, then $\eta_x$ is an isometric bijection. 
\eit 
\item $d_S^X \leq \nu_H(S, T) \,d_T^X$ \ \ $(S, T \in {\cal S}(G))$ \hspp (Here, we set $\infty \cdot 0 = 0$.) 
\eenum 
\efact 

Recall the bijections ${\cal S}'(G) \approx {\cal S}''(G) \approx \widetilde{\cal S}(G)$ in Notation~\ref{not_S_p}\,[II](4). 

\bpropb 
Consider the map 
$\bary[t]{c@{ \ }c@{ \ }c@{ \ }c@{ \ }c@{ \ }c@{}}
\phi : & ({\cal S}''(G)^\ast,  \widehat{\rho}\,) & \ \approx \ & 
(\widetilde{\cal S}(G)^\ast,  \widetilde{\widehat{\rho}\,} = \widehat{\widetilde{\rho}\,} ) & \lra & 
(F(X^2 - \Delta(X),\widetilde{\IR}_{\geq 0}), \widehat{\mu}) \\[0.5mm] 
& S & & [S] & & \log \,d_S^X|
\eary$. \\[1mm] 
Here, $d_S^X| := d_S^X|_{X^2 - \Delta(X)}$. Suppose $X$ is not a singleton. 
\benum 
\item $\phi$ is 1-Lipschitz. 
\item If $\eta_x : G \to X$ is injective for some $x \in X$, then 
\bit 
\itemI $\lambda(d_T^X|, d_S^X|) = \nu_H(S,T)$ \ \ $(S, T \in {\cal S}(G)^\ast)$. 
\itemII $\phi$ is an isometric injection. ($\phi$ is an anti-isometric injection with respect to $\rho$ and $\mu$). 
\eit 
\eenum 
\eprop

\bpfb Let $r := \lambda(d_T^X|, d_S^X|) \in \widetilde{\IR}_{\geq 1}$. 
\benum 
\item $r \leq \nu_H(S, T)$ \ \ by Fact~\ref{fact_eta}\,(2). 
\item 
\bit 
\itemI $d_S^X| \leq r \,d_T^X|$ \hsh \tf \ $d_S^G| \leq r \,d_T^G|$ \ \ by Fact~\ref{fact_eta}\,(1)(ii). 
\hsh \tf \ $\nu_S| \leq r\,\nu_T|$ \\
\tf \ $\nu_H(S, T) = \sup \, \nu_S(T) \leq r \,\sup\,\nu_T(T) = r$
\eit 
\eenum 
\vskip -7mm 
\epf 

\subsection{Symmetry of $(X, G, \rho)$ } \mbox{} 

\bnotb 
\benum 
\item Suppose $(X, G, \rho)$ and $(Y, H, \zeta)$ are right group actions on sets $X$ and $Y$. \hsp 
\smash{\raisebox{0mm}{$\xymatrix@M+1pt{
X \times G \ar[d]_-{\rho} \ar[r]^-{f \times \phi} & Y \times H\ar[d]^-{\zeta} \\
X \ar[r]^-{f} & Y
}$}} \\
An equivariant map \ $(f, \phi) : (X, G, \rho) \lra (Y, H, \zeta)$ \ is \\
\btab[t]{@{}l}
a pair of a map $f : X \to Y$ and a group homomorphism $\phi : G \to H$ \\[1.5mm]  %\ or \ an anti gr. homo. \\[1.5mm] 
\hsp such that \ $f(xa) = f(x)\phi(a)$ \ $(\upfa (x,a) \in X \times G)$.  
\etab \\[2mm] 
The composition of equivariant maps and the identity maps are defined as follows. 
\vskip 1mm 
\bit 
\itemI $(g, \psi) (f, \phi) := (gf, \psi \phi)$ \ for \ 
$\bary[b]{@{}c@{ \ }c@{ \ }c@{ \ }c@{ \ }c@{}}
& \mbox{\small $(f, \phi)$} & & \mbox{\small $(g, \psi)$} & \\[-1mm] 
(X, G, \rho) & \lra & (Y, H, \zeta) & \lra & (Z, I, \eta) 
\eary$. 
\itemII $\id_{(X, G, \rho)} := (\id_X, \id_G)$
\eit 
\item In the category of right group actions on sets and equivariant maps between them.  
\bit 
\itemI ${\rm End}(X, G, \rho) \supset {\rm Aut}(X, G, \rho)$ denote the monoid of endomorphisms and the group of automorphisms of 
$(X, G, \rho)$. 

\itemII An equivariant map $(f, \phi) : (X, G, \rho) \lra (Y, H, \zeta)$ is an isomorphism if and only if both $f$ and $\phi$ are bijective. 
In this case $(f,\phi)^{-1} = (f^{-1},\phi^{-1})$. 
\eit 
\eenum 
\enot 

Suppose $(X, G, \rho)$ is a right group action on a set $X$, 
$M$ is a set with a binary operation and \\ 
$(f, \phi) : M \ni u \lmt (f_u, \phi_u) \in {\rm End}(X, G, \rho)$ is an anti-homomorphism. 

\bnotb Let $u,v \in M$, $x \in X$, $a,b \in G$ \ and \ $S \subset G$. 
\benum 
\item $\phi : M \lra {\rm End}(G)$ is an anti-homomorphism. 
\hsp $\circ$ \ $(f_{uv}, \phi_{uv}) = (f_v, \phi_v)(f_u, \phi_u) = (f_vf_u, \phi_v\phi_u)$ 
\item $x^u := f_u(x)$, \ $a^u := \phi_u(a)$ 
\bit 
\itemI $x^{uv} = (x^u)^v$ \hsp 
(ii) \ $(xa)^u = x^u a^u$ \hsp 
\itemiii $a^{uv} = (a^u)^v$ \hsp 
(iv) \ $(ab)^u = a^u b^u$ \hsp 
(v) \ $(a^b)^u = (a^u)^{b^u}$ 
\eit 
\item \bit 
\itemI $S^u := \phi_u(S) = \{ a^u \mid a \in S \}$ \hsp (ii) \ 
$S \subset G$ is $\phi$-invariant. \LLRA $S^u \subset S$ \ $(\upfa u\in M)$ 
\eit 
\eenum 
\enot

\bfactb Let $A \subset X$, $S, T \subset G$, $u,v \in M$ and $n \in \IZ_{\geq 0}$. 
\benum 
\item (i) \ $(A^u)^v = A^{uv}$ \hsp (ii) \ $(AS)^u = A^uS^u$  
\item (i) \ $(S^u)^v = S^{uv}$ \hsp \,(ii) \ $(ST)^u = S^uT^u$ \hsp 
(iii) \ $(S^T)^u = (S^u)^{T^u}$ 
\item[(3)] (i) \ $(S^n)^u = (S^u)^n$ \hsp 
(ii) \ $(S^{\leq n})^u = (S^u)^{\leq n}$ \hsp  
(iii) \ $\big( S^\infty \big)^u = (S^u)^\infty$ 
\eenum 
\efact 

\bfact\label{fact_sym_d_S} Let $S \in {\cal S}(G)$. 
\benum  
\item $d_{S^u}^X(x^u, y^u) \leq d_S^X(x, y)$ \ \ $(\upfa x,y \in X, \ u \in M)$
\item If $S$ is $\phi$-invariant, then \ $d_S^X(x^u, y^u) \leq d_S^X(x, y)$ \ $(\upfa x,y \in X, \ u \in M)$.  
\item Suppose $M$ is a monoid, $e \in M$ is the unit element and 
$(f,\phi)$ is an anti-homomorphism of monoids \ (i.e., $(f_e, \phi_e) = \id$ \ \ $(x^e = x, a^e = a)$). 
Let $U(M)$ denote the group of units of $M$. For $u \in U(M)$ the following holds. 
\bit 
\itemI $d_{S^u}^X(x^u, y^u) = d_S^X(x,y)$ \ \ $(\upfa x,y \in X)$
\itemII If $S$ is $\phi$-invariant, then $S^u = S$ and $d_S^X(x^u, y^u) = d_S^X(x,y)$ \ $(\upfa x,y \in X)$. 
\eit 
\eenum 
\efact 

\bexp\label{exp_sym} \mbox{} 
\benum 
\item Suppose $M$ is a $\ast$-subset of ${\rm End}(X, G, \rho)$ (i.e., $M$ is closed under the multiplication). 
We reverse the multiplication of $M$ so that the inclusion $M \subset {\rm End}(X, G, \rho)$ is an anti-homomorphism. 
\bit 
\itemI $M$ acts on $X$ and $G$ by \ $x^{(f,\phi)} = f(x)$ \ and \ $a^{(f, \phi)} = \phi(a)$ \ $(x \in X, a \in G, (f, \phi) \in M)$.   
\itemII If $S \subset G$ is $M$-invariant (i.e., $S^{(f, \phi)} = \phi(S) \subset S$ \ $(\upfa (f, \phi) \in M)$), then \\
\hsp $d_S^X(f(x), f(y)) \leq d_S^X(x, y)$ \ \ $(\upfa x,y \in X, \ (f, \phi) \in M)$.  
\eit 
\item When $M \subset {\rm Aut}(X, G, \rho)$ is a subgroup under the reverse multiplication, 
\bit 
\itemI $d_{\phi(S)}^X(f(x), f(y)) = d_S^X(x, y)$ \ \ $(\upfa x,y \in X, \ (f, \phi) \in M)$ 
\itemII If $S \subset G$ is $M$-invariant, then $\phi(S) = S$ and \\
\hsp $d_S^X(f(x), f(y)) = d_S^X(x, y)$ \ \ $(\upfa x,y \in X, \ (f, \phi) \in M)$.  
\eit 
\eenum
\eexp 

\bexpb Each action $(X, G, \rho)$ admits a canonical anti-homomorphism \\ 
\hsppp $(f, \phi) : G \lra {\rm Aut}(X, G, \rho)$ : $(f_a, \phi_a) = (\rho_a, a^{-1}( \ )a)$. \\
However, this symmetry yields no new information, as seen below. 
\benum
\item Let $x \in X, \ a, g \in G$ \ and \ $S \subset G$. 
\bit 
\itemI $x^a = \rho_a(x) = xa$, \ \ $g^a = \phi_a(g) = a^{-1}ga = g^a$ \ (in the sense of Section 3) 
\itemII $S^a = \phi_a(S) = a^{-1}Sa = S^a$ (in the sense of Section 3)
\itemiii $S$ is $\phi$-invariant. \LLRA $a^{-1}Sa \subset S$ \ $(\upfa a \in G)$ \LLRA $S$ is conjugation invariant.  
\eit 
\item  
Fact~\ref{fact_sym_d_S}\,(3) implies the following statements, which exactly restate Fact~\ref{fact_d_S}\,(3). 
\bit 
\itemI $d_{S^a}^X(xa,ya) = d_S^X(x,y)$ \ \ $(\upfa x,y \in X, \ \upfa a \in G)$.  
\itemII If $S$ is conjugation invariant, then \ $d_S^X(xa,ya) = d_S^X(x,y)$ \ $(\upfa x,y \in X, \ \upfa a \in G)$. 
\eit 

\item If we consider the action of a group $G$ on itself by the right translation \\
\hsp \hsh $(G, G, R)$ : \ \ $R : G \caleft G$ : $R(x,a) = xa$ \\
and the anti group homomorphism \\
\hsp \hsh $(f, \phi) : G \lra {\rm Aut}(G, G, R)$ : $(f_a, \phi_a) = (R_a, a^{-1}(\ )a)$, \\
then we recover the statements in Fact~\ref{fact_d-nu}.
\eenum  
\eexp 

\section{Word length on $\ast$-sets} 

\subsection{$\ast$-sets} \mbox{} 

In this section we extend the arguments for groups in the previous sections to the case of sets with binary operations.

\bdefnb A $\ast$-set means a set $X$ with binary operations $\ast_i : X \times X \to X$ \ $(i= 1, \cdots, m)$.   
\edefn 

Let $X = (X, (\ast_i)_{i \in [m]})$ be a nonempty $\ast$-set.  

\bnotb 
\benum
\item A $\ast$-subset of $X$ means a subset $Y \subset X$ such that $y_1 \ast_i y_2\in Y$ \ $(\upfa i \in [m])$ for any $y_1, y_2 \in Y$.
\bit 
\itemI In the case of unital $\ast$-sets (i.e., $\ast$-sets with an identity element),  we can consider unital $\ast$-subsets. 
\eit 
\item ${\cal S}_\ast(X) : =$ the set of $\ast$-subsets of $X$ 
\bit 
\itemI ${\cal S}_\ast(X)$ is closed under intersection 
(i.e., $Y_\lambda \in {\cal S}_\ast(X)$ \ $(\lambda \in \Lambda)$ \LRA 
$\bigcap_{\lambda \in \Lambda} Y_\lambda \in {\cal S}_\ast(X)$). 
\eit 

\item Let $x,y, x_1, \cdots, x_n \in X$. 
\bit 
\itemI The symbol $x \ast y$ denotes any of $x \ast_i y$ \ $(i \in [m])$. 
\itemII The symbol $x_1 \ast \cdots \ast x_n$ denotes any product of $x_1, \cdots, x_n$ in this order 
under the (mixed) operations $\ast_i$ $(i\in [m])$ and suitable $( \ )$'s. 
\itemC ${\cal S}_\ast(X) \ni Y \ni x_1, \cdots, x_n$ $(n \in \IZ_{\geq 1})$ \LRA $x_1 \ast \cdots \ast x_n \in Y$ 
\eit 
\item $A \ast B := \{ a \ast b \mid a \in A, b \in B \} \equiv \{ a \ast_i b \mid a \in A, b \in B, i \in [m] \}$ \hsh $(A, B \subset X)$ 
\item[] $A_1 \ast \cdots \ast A_n := \{ a_1 \ast \cdots \ast a_n \mid a_i \in A_i \ (i \in [n]) \}$ \hsh 
$(n \in \IZ_{\geq 1}, \ A_1, \cdots, A_n \subset X)$ 
\eenum 
\enot

\bnotb Let $S \subset X$. 
\benum
\item The $\ast$-subset of $X$ generated by $S$ means the smallest $\ast$-subset of $X$ including $S$. 
It exists and is obtained formally by $\langle S \rangle := \bigcap \{ Y \in {\cal S}_\ast(X) \mid S \subset Y \}$.  
\item The $\ast$-subset $\langle S \rangle$ has the following explicite description : \ $\langle S \rangle = S^\infty$ 
\bit 
\itemI $S^n := \underbrace{S \ast \cdots \ast S}_{n} 
\equiv \{ x_1 \ast \cdots \ast x_n \mid x_i \in S \ (i \in [n]) \}$ \ \ $(n \in \IZ_{\geq 1})$
\vskip 2mm 
\item[] $\circ$ \ Inductive definition : \\
\hsp $S^1 = S$ \hsp $S^k$ \ $( k\in [n])$ \ \ $\rsa$ \ \ $S^{n+1} = \bigcup \{ S^k \ast S^\ell \mid k, \ell \in \IZ_{\geq 1}, k+\ell = n+1 \}$ 
\item[] $\circ$ \ $S^2 = S \ast S$, \ \ $S^3 = (S^1 \ast S^2) \cup (S^2 \ast S^1) = (S \ast (S \ast S)) \cup ((S \ast S) \ast S)$, \ etc.

\itemII $S^{\leq n} := \bigcup_{k \in [n]} S^k$ \ \ $(n \in \IZ_{\geq 1})$ \hsp $S^\infty := \bigcup_{n \in \IZ_{\geq 1}} S^n$ 
\vskip 0.5mm
\itemiii In the unital case (i.e., $X$ has the identity element $e$), we set $S^0 := \{ e \}$ and $S^\infty := \bigcup_{n \in \IZ_{\geq 0}} S^n$.
\eit 
\eenum 
\enot

\bfactb 
\benum 
\item \bit 
\itemI $S^k \ast S^\ell \subset S^{k + \ell}$ \ \ $(k, \ell \in \IZ_{\geq 1})$ 
\hspp $\circ$ \ This implies that $S^\infty \in {\cal S}_\ast(X)$ and $\langle S \rangle = S^\infty$. 
\itemII $S^{\leq k} \ast S^{\leq \ell} \subset S^{\leq k+\ell}$ \ \ $(k, \ell \in \IZ_{\geq 1})$ 
\eit 
\item  
\bit 
\itemI $(S^k)^\ell \subset S^{k\ell}$ \hsh (ii) \ $(S^{\leq k})^{\leq \ell} \subset S^{\leq k\ell}$ \hsp $(k, \ell \in \IZ_{\geq 1})$
\eit 
\eenum 
\efact 

\subsection{Word length on a $\ast$-set} \mbox{} 

\bdefnb (Word length) \ Let $S \subset X$.  
\benum 
\item $\nu_S : X \lra \widetilde{\IZ}_{\geq 1}$ : \ 
$\nu_S(x) = 
\left\{ \hspace*{-1mm} \bary[c]{ll}
\min \{ n \in \IZ_{\geq 1} \mid x \in S^n \} & (x \in S^\infty) \\[1mm]
\infty & (x \in G - S^\infty)  
\eary \right.$ 
\vskip 2mm 
\bit 
\itemI In the unital case $\IZ_{\geq 1}$ is replaced by $\IZ_{\geq 0}$.
\eit 
\eenum 
\edefn 

\bfactb
\benum 
\item 
\bit 
\itemI $\nu_S(x) \leq n$ \LLRA $x \in S^{\leq n}$ \ \ $(x \in X, n \in \IZ_{\geq 1})$
\itemII $\nu_S^{-1}(1) = S$ \hsh 
$(\nu_S^{-1}(1) = S - \{ e \}$ in the unital case)
\eit 
\item 
\bit 
\itemI $\nu_S(x \ast y) \leq \nu_S(x) + \nu_S(y)$ \ \ $(x,y \in X)$ 
\itemII $\nu_S(x_1 \ast \cdots \ast x_n) \leq \nu_S(x_1) + \cdots + \nu_S(x_n)$ \hsh $(x_i \in X \ (i \in [n]))$ \hsh  
$(n \in \IZ_{\geq 1})$ 
\eit 
\eenum
\efact 

\bdefnb (the induced metric on ${\cal S}(X)^\times = {\cal S}(X) - \{ \emptyset \}$)  
\benum
\item $\nu_H : {\cal S}(X)^\times \times {\cal S}(X)^\times \lra \widetilde{\IZ}_{\geq 1}$ : \ $\nu_H(S, S') =  \sup \nu_S(S')$ 
\item $\rho := \log \nu_H : {\cal S}(X)^\times \times {\cal S}(X)^\times \lra \log \widetilde{\IZ}_{\geq 1} \subset \widetilde{\IR}_{\geq 0}$ 
\eenum 
\edefn

\bfactb 
\benum 
\item $\nu_H(S, S') \leq n$ \LLRA $S' \subset S^{\leq n}$ \ \ $(S, S' \in {\cal S}(X)^\times, n \in \IZ_{\geq 1})$ 
\item $\nu_H \in {\cal M}_{m,as}({\cal S}(X)^\times)$, \ \ $\nu_H$ is nondegenerate. 
\bit 
\itemI $\nu_H(S, S') = 1$ \LLRA $S' \subset S$ \hsp (i)$'$ \ $\nu_H(S, S') = \nu_H(S', S) = 1$ \LLRA $S = S'$
\itemII $\nu_H(S, S'') \leq \nu_H(S, S')\nu_H(S',S'')$ \ \ $(S,S',S'' \in {\cal S}(X)^\times)$ 
\eit 
\item $\rho \in {\cal M}_{as}({\cal S}(X)^\times)$, \ \ $\rho$ is nondegenerate. \hsp \tf $\widehat{\rho} \in {\cal M}({\cal S}(X)^\times)$
\eenum 
\efact 

\bexp\label{exp_*-set_01} 
Some examples of subspaces ${\cal X} \subset {\cal S}(X)^\times$ with $\widehat{\rho}|_{{\cal X}} < \infty$ \\ 
\hsp Let ${\cal S}_g(X) := \{ S \in {\cal S}(X) \mid X = S^\infty \} \subset {\cal S}(X)^\times$.  
\benum 
\item ${\cal S}_{fg}(X) := \{ S \in {\cal S}_g(X) \mid |S| < \infty \}$
\item ${\cal S}_g(X; \, \mbox{finite  rel $S_0$}) := \{ S \in {\cal S}_g(X) \mid |S - S_0| < \infty, |S_0 - S| < \infty \}$ \hsh 
$(S_0 \in {\cal S}_g(X))$ 

\bit 
\itemI If $S \in {\cal S}_g(X)^\times$, $S' \in {\cal S}(X)^\times$ and $|S' - S| < \infty$, then $\nu_H(S, S') < \infty$ 
\eit 
\eenum 
\eexp

For word metrics on the $\ast$-set $X$ we have the following candidates. 

\bnotb (Word metrics) \ Let $S \subset X$. \\[2mm]
\hsp \hsh  
$d_S : X \times X \lra \widetilde{\IZ}_{\geq 0}$ : 
$d_S(x,y) = 
\left\{ \hspace*{-1mm} \bary[c]{ll}
\min \{ n \in \IZ_{\geq 0} \mid y \in [x, S, n] \} \ & (y \in 
\bigcup_{n \in \IZ_{\geq 0}} [x, S, n]) \\[2mm]
\infty & (\mbox{otherwise})  
\eary \right.$ 
\vskip 3mm 
\hsh Here, we can take the family $[x, S, n]$ \ $(x \in X, n \in \IZ_{\geq 0})$ as follows : 
\benum
\item[] 
\bit 
\item[(1)] $[x, S, n] := x \ast \underbrace{S \ast \cdots \ast S}_{n}$ \ \ \ $(n \in \IZ_{\geq 1})$ \hsp $[x,S,0] := \{ x \}$ 
\hsp $\circ$ \ $[x, S, n] \ \supset \ x \ast (S^n)$ 
\vskip 2mm 
\item[(2)] $[x, S, n]' := ((\cdots ((x \ast \underbrace{S) \ast S) \ast \cdots \ast S )\ast S}_{n}$ \ \ \ $(n \in \IZ_{\geq 1})$ \hsp $[x,S,0]' := \{ x \}$ 
\eit 
\eenum
\vskip 3mm 
\hsp We denote by $d_S'$ the metric $d_S$ defined in the case (2). 
\enot 

\bfactb $d_S, d_S' \in {\cal M}_{as}(X)$ and they are nondegenerate. \hsh \tf \ $\widehat{d_S}, \,\widehat{d_S'} \in {\cal M}(X)$. 
\benum 
\item 
\bit 
\itemI $d_S(x,y) = 0$ \LLRA $x = y$ \hsh (ii) \ $d_S(x,z) \leq d_S(x,y) + d_S(y,z)$ \hsh $(\upfa x,y,z \in X)$ 
\eit 
\item 
\bit 
\itemI $d_S'(x,y) = 0$ \LLRA $x = y$ \hsh (ii) \ $d_S(x,z) \leq d_S(x,y) + d_S(y,z)$ \hsh $(\upfa x,y,z \in X)$
\eit 
\eenum
\efact

\subsection{Symmetry of a $\ast$-set} \mbox{} 

Let $X$ be a $\ast$-set. 

\bnotb  
\benum
\item A map $f : X \to X$ is said to be a $\ast$-endomorphism if $f(x \ast_i y) = f(x) \ast_i f(y)$ \ $(x,y \in X, \ i \in [m])$).
\item[] $\circ$ \ In this case, \ (i) \ $S \in {\cal S}_\ast(X)$ \LRA $f(S) \in {\cal S}_\ast(Y)$ \hsh 
(ii) \ $R \in {\cal S}_\ast(Y)$ \LRA $f^{-1}(R) \in {\cal S}_\ast(X)$
\item ${\rm End}(X)$ denotes the monoid of $\ast$-endmorphisms of $X$
\eenum 
\enot 

Suppose $(R, \star)$ is a nonempty set with a binary operation and $\phi : R \lra {\rm End}(X)$ is an anti-homomorphism. 

\bnotb 
 \benum 
\item $x^a := \phi(a)(x) \in X$ \ \ $(x \in X, a \in R)$ \hsh 
(2) \ $S^A := \{ x^a \mid x \in S, a \in A \}$ \ \ $(S \subset X, A \subset R)$ 
 \item $S \subset X$ is $\phi$-invariant \LLRAdefn $S^a \subset S$ \ \ $(\upfa a \in R)$ \LLRA $S^R \subset S$   
\eenum
\enot 

\bfactb \ Let $x,y \in X$, $S \subset X$, $a,b \in R$ and $A, B \subset R$. 
\benum 
 \item 
 \bit 
 \itemI $(x \ast y)^a = x^a \ast y^a$ \hsh (ii) \ $(x^a)^b = x^{a \star b}$ \hsh (iii) \ $(S^A)^B = S^{A \star B}$ 
 \hsh (iv) \ $(S \cup T)^A = S^A \cup T^A$ 
 \eit 
\vskip 0.5mm 
\item 
\bit 
\itemI $(Y \ast Z)^a = Y^a \ast Z^a$ \ \ $(Y, Z \subset X)$ \hsp (ii) \ 
$(\bigcup_{\lambda \in \Lambda} Y_\lambda)^a = \bigcup_{\lambda \in \Lambda} Y_\lambda^a$ \ \ 
$(Y_\lambda \subset X \ (\lambda \in \Lambda))$
\eit 
\vskip 1mm 
\item 
\bit 
\itemI $(x_1 \ast \cdots \ast x_n)^a = x_1^a \ast \cdots \ast x_n^a$ \ \ $(x_i \in X \ (i \in [n]))$ \\
 \ \ (The products are taken in the same form in the left side and right side.)  
\itemII $(S^n)^a = (S^a)^n$ \ \ $(n \in \IZ_{\geq 1})$ \hsp $(S^\infty)^a = (S^a)^\infty$ \hsp $\langle S \rangle^a = \langle S^a \rangle$
\eit 
\item 
\bit 
\itemI $S^R$ is $\phi$-invariant. 
\hsh (ii) \ $(S \cup S^R)^R = S^R \subset S \cup S^R$ \hsh \tf \ $S \cup S^R$ is $\phi$-invariant.

\itemII $S$ is  $\phi$-invariant. \LRA $(S^n)^a \subset S^n$ \ $(n \in \IZ_{\geq 1})$, \ $(S^\infty)^a \subset S^\infty$
\eit 
\eenum
\efact

\bfactb Suppose $S \in {\cal S}(X)^\times$ is $\phi$-invariant. 
\benum 
\item $\nu_S(x^a) \leq \nu_S(x)$ \ \ $(\upfa x \in X, \upfa a \in R)$ 
\item $\nu_H(S, T^R) \leq \nu_H(S, T \cup T^R) = \nu_H(S, T)$ \ \ $(\upfa T \in {\cal S}(X)^\times)$ 
\eenum 
\efact 

\bexpb Further examples of subspaces ${\cal X} \subset {\cal S}(X)^\times$ with $\widehat{\rho}|_{{\cal X}} < \infty$  
\benum
\item ${\cal S}_{\phi\,\mbox{-}fg}(X) 
:= \{ S \in {\cal S}_g(X) \mid S = F^R \ {\rm or} \ F \cup F^R \ \ (\upexi F \in {\cal S}_f(X)^\times) \}$ 

\item[] $\circ$ \ For any $S, T \in {\cal S}_{\phi\,\mbox{-}fg}(X)$
\bit 
\itemI $S \in {\cal S}_g(X)$, $S$ is $\phi$-invariant. \hsh (ii) \ $T = F^R$ or $F \cup F^R$ for some $F \in {\cal S}_f(X)^\times$. 
\itemiii $\nu_H(S, F^R) \leq \nu_H(S, F \cup F^R) = \nu_H(S, F) < \infty$ \hsh \tf \ $\nu_H(S, T) < \infty$ 
\eit 
\vskip 1mm 

\item ${\cal S}_{\phi\,\mbox{-}g}(X; \, \mbox{finite  rel $F_0$}) 
:= \{ S \in {\cal S}_g(X) \mid S = F^R \ {\rm or} \ F \cup F^R \ \ (\upexi F \in {\cal S}(X, \, \mbox{finite  rel $F_0$})^\times \}$ 

\item[] $\circ$ \ ${\cal S}(X; \, \mbox{finite  rel $F_0$})^\times := \{ F \in {\cal S}(X)^\times \mid |F - F_0| < \infty, |F_0 - F| < \infty \}$ 
\hsp $(F_0 \in {\cal S}(X)^\times)$
\item[] $\circ$ \ For any $S, T \in {\cal S}_{\phi\,\mbox{-}g}(X; \, \mbox{finite  rel $F_0$})$ 
\bit 
\itemI $S \in {\cal S}_g(X)$, $S$ is $\phi$-invariant. \hsf 
(ii) \ $S = F \cup F^R$ for some $F \in {\cal S}(X, \, \mbox{finite  rel $F_0$})^\times$. 
\itemiii $T = G \cup G^R$ for some $G \in {\cal S}(X, \, \mbox{finite  rel $F_0$})^\times$. 
\itemiv $G - S \subset G  - F \subset (G - F_0) \cup (F_0 - F) \in {\cal S}_f(X)$ \\
\tf \ $\nu_H(S, G) < \infty$ \hsh (cf. Example~\ref{exp_*-set_01}\,(2)(i)) 
\itemv $\nu_H(S, G^R) \leq \nu_H(S, G \cup G^R) = \nu_H(S, G) < \infty$ \hsh \tf \ $\nu_H(S, T) < \infty$ 
\eit 
 \eenum 
 \eexp

\edc

%% file: macro.tex
\theoremstyle{definition}

\newtheorem{thm}{Theorem}[section]
\newtheorem{prop}{Proposition}[section]
\newtheorem{cor}{Corollary}[section]
\newtheorem{lem}{Lemma}[section]

\theoremstyle{definition}

\newtheorem{remark}{Remark}[section]
\newtheorem{claim}{Claim}[section]
\newtheorem{defn}{Definition}[section]
\newtheorem{notation}{Notation}[section]
\newtheorem{example}{Example}[section]
\newtheorem*{problem}{Problem}%[section] 
\newtheorem{fact}{Fact}[section]
\newtheorem*{property}{Property}%[section]
\newtheorem*{observation}{Observation}%[section]
\newtheorem*{assumption_*}{Assumption $(\ast)$} 

\newtheorem*{fact and notation I}{Fact and Notation I}
\newtheorem*{fact and notation II}{Fact and Notation II}

\def \bthm {\begin{thm}}
\def \bthmb {\begin{thm} \mbox{}}
\def \ethm {\end{thm}} 

\def \bprop {\begin{prop}}
\def \bpropb {\begin{prop} \mbox{}}
\def \eprop {\end{prop}} 

\def \blem {\begin{lem}}
\def \blemb {\begin{lem} \mbox{}}
\def \elem {\end{lem}} 

\def \bcor {\begin{cor}}
\def \bcorb {\begin{cor} \mbox{}}
\def \ecor {\end{cor}}

\def \bdefn {\begin{defn}} 
\def \bdefnb {\begin{defn} \mbox{}}
\def \edefn {\end{defn}} 

\def \bfact {\begin{fact}} 
\def \bfactb {\begin{fact} \mbox{}}
\def \efact {\end{fact}} 

\def \bnot {\begin{notation}} 
\def \bnotb {\begin{notation} \mbox{}}
\def \enot {\end{notation}} 

\def \brem {\begin{remark}} 
\def \bremb {\begin{remark} \mbox{}} 
\def \erem {\end{remark}} 

\def \bpty {\begin{property}} 
\def \bptyb {\begin{property} \mbox{}} 
\def \epty {\end{property}} 

\def \bexp {\begin{example}} 
\def \bexpb {\begin{example} \mbox{}} 
\def \eexp {\end{example}} 

\def \bexp {\begin{example}} 
\def \bexpb {\begin{example} \mbox{}} 
\def \eexp {\end{example}} 

\def \bcl {\begin{claim}}
\def \bclb {\begin{claim} \mbox{}} 
\def \ecl {\end{claim}} 

\def \bobs {\begin{observation}}
\def \bobsb {\begin{observation} \mbox{}} 
\def \eobs {\end{observation}} 

\def \btab {\begin{tabular}}
\def \btabb {\begin{tabular} \mbox{}} 
\def \etab {\end{tabular}} 

\def \bary {\begin{array}}
\def \baryb {\begin{array} \mbox{}} 
\def \eary {\end{array}} 

\def \barr {\begin{array}}
\def \barrb {\begin{array} \mbox{}} 
\def \earr {\end{array}} 

\def \bpf {\begin{proof}}
\def \bpfb {\begin{proof} \mbox{}} 
\def \epf {\end{proof}} 

\def \benum {\begin{enumerate}}
\def \eenum {\end{enumerate}} 

\def \bit {\begin{itemize}}
\def \eit {\end{itemize}} 

\def \balign {\begin{align*}}
\def \ealign {\end{align*}}

\def \bCD {\begin{CD}}
\def \eCD {\end{CD}}

\def \edc {\end{document}} 

\def \beq {\begin{equation}}
\def \eeq {\end{equation}}
\def \itemC {\item[$\circ$]}

\newcommand{\bmp}[1]{\begin{minipage}[t]{#1} \baselineskip 6mm}
\def \emp {\end{minipage}}

\def \bec {\mbox{($\because$) }}
\def \tf {\mbox{$\therefore$ }}

\def \upexi {\mbox{${}^\exists$}} 

\def \upfa {\mbox{${}^\forall$}}

\def \lmt {\longmapsto}

\newcommand{\relsb}[3]{
\bary[c]{@{}c@{}}
\mbox{\scriptsize $#1$} \\[-1.5mm] 
\mbox{$#2$} \\[-1.5mm] 
\mbox{\scriptsize $#3$}
\eary
}

\def \cal {\mathcal}
\def \phi {\varphi}

\renewcommand{\Theta}{\varTheta}
\renewcommand{\Sigma}{\varSigma}
\renewcommand{\Phi}{\varPhi}
\renewcommand{\Psi}{\varPsi}
\renewcommand{\Lambda}{\varLambda}
\renewcommand{\Gamma}{\varGamma}
\renewcommand{\rho}{\varrho}

\def \Llra {\Longleftrightarrow}
\def \Lra {\Longrightarrow}
\def \Lla {\Longleftarrow}

\def \lra {\longrightarrow}
\def \lla {\longleftarrow}
\def \LLRA {\ $\Llra$ \ } 
\def \LRA {\ $\Lra$ \ }

\def \rsa {\rightsquigarrow}

\def \caleft {\curvearrowleft}

\def \LLRAdefn
{\btab[b]{c}
\mbox{\footnotesize def} \\[-1.5mm] 
\ $\Llra$ \ 
\etab}

\def \LLRAequiv
{\btab[b]{c}
\mbox{\footnotesize equiv} \\[-1.5mm] 
\ $\Llra$ \ 
\etab}

\def \itemi {\item[(i)\,\mbox{}]}
\def \itemii {\item[(ii)]}
\def \itemiii {\item[(iii)]}
\def \itemiv {\item[(iv)]}
\def \itemv {\item[(v)\,\mbox{}]}
\def \itemI {\item[(i)\ \mbox{}]}
\def \itemII {\item[(ii)\,\mbox{}]}

\def \itema {\item[(a)]} 
\def \itemb {\item[(b)]}

\def \hsf {\hspace*{2mm}}
\def \hsh {\hspace*{5mm}}
\def \hsp {\hspace*{10mm}}
\def \hspp {\hspace*{20mm}}
\def \hsppp {\hspace*{30mm}}

\def \bprob {\begin{problem}}
\def \bprobb {\begin{problem} \mbox{}}
\def \eprob {\end{problem}}

\def \IR {{\Bbb R}}

\def \IZ {{\Bbb Z}}

\newcommand{\id}{\mathrm{id}}

\newcommand{\incto}{\ar@{^{(}->}}
\newcommand{\into}{\ar@{>->}}
\newcommand{\onto}{\ar@{->>}}